\documentclass[11pt]{amsart}

\usepackage{amsmath,amssymb,latexsym,amsthm,mathdots,yhmath}

\usepackage{graphicx,amsfonts,subfigure}

\newtheorem{theorem}{Theorem}[section]
\newtheorem{lemma}[theorem]{Lemma}
\newtheorem{proposition}[theorem]{Proposition}

\theoremstyle{remark}
{
    \newtheorem{definition}[theorem]{Definition}

\newtheorem{assumption}[theorem]{Assumption}
}

\newtheorem{remarks}{Remark}[section]{\bf}{\rmfamily}

\def\Ker{\mathop{\rm Ker}}

\newcommand{\shoot}{\mathcal{S}}

\newcommand{\tras}{^\top}

\newcommand{\lam}{[\lambda]}
\newcommand{\lamh}{[\hat\lambda]}
\newcommand{\lamLQ}{[\lambda^{LQ}]}

\newcommand{\Uspace}{\U }

\newcommand{\intT}{\int_0^T }
\newcommand{\gr}{>}
\newcommand{\mi}{<}

\def\half{\mbox{$\frac{1}{2}$}}

\newcommand{\ddt}{\frac{\rm d}{{\rm d} t} }

\def\ph{\hat{p}}

\def\uh{\hat{u}}

\def\wh{\hat{w}}
\def\xh{\hat{x}}

\def\Ih{\hat{I}}
\def\Th{\hat{T}}
\def\Wh{\hat{W}}

\def\B{\mathcal{B}}
\def\C{\mathcal{C}}

\def\F{\mathcal{F}}
\def\G{\mathcal{G}}
\def\H{\mathcal{H}}
\def\I{\mathcal{I}}
\def\J{\mathcal{J}}
\def\K{\mathcal{K}}
\def\L{\mathcal{L}}
\def\M{\mathcal{M}}

\def\P{\mathcal{P}}

\def\T{\mathcal{T}}
\def\U{\mathcal{U}}

\def\W{\mathcal{W}}
\def\X{\mathcal{X}}
\def\Y{\mathcal{Y}}

\newcommand{\cR}{I\!\! R}

\newcommand\be{\begin{equation}}
\newcommand\ee{\end{equation}}
\newcommand{\benl}{\begin{equation*}}
\newcommand{\eenl}{\end{equation*}}
\newcommand\ba{\begin{array}}
\newcommand\ea{\end{array}}
\newcommand{\bean}{\begin{eqnarray*}}
\newcommand{\eean}{\end{eqnarray*}}

\def\ds{\displaystyle}

\newcommand{\dtt}{\mathrm{d}t}


\title[A Shooting Algorithm for Problems with Singular Arcs]{A Shooting Algorithm for Optimal Control Problems with Singular Arcs$^{1,2}$}

\author[M.S. Aronna]{M. Soledad Aronna}
\address{M.S. Aronna\\ Fellowship within the ITN Marie Curie Network SADCO at Universit\`a degli Studi di Padova\\  Padova  35121, Italy}
\email{aronna@math.unipd.it}

\author[J.F. Bonnans]{J. Fr\'ed\'eric Bonnans}
\address{J.F. Bonnans\\ INRIA-Saclay and CMAP,  \'Ecole~Polytechnique, 91128~Palaiseau, France}
\email{Frederic.Bonnans@inria.fr}

\author[P. Martinon]{Pierre Martinon}
\address{P. Martinon\\ INRIA-Saclay and CMAP,  \'Ecole~Polytechnique, 91128~Palaiseau, France}
\email{Pierre.Martinon@inria.fr}

\begin{document}

\maketitle

\footnotetext[1]{This work is supported by the European Union under the 7th Framework Programme «FP7-PEOPLE-2010-ITN»  Grant agreement number 264735-SADCO}

\footnotetext[2]{This article was accepted for publication in {\it Journal of Optimization, Theory and Applications.}}

\begin{abstract}
In this article, we propose a shooting algorithm for a class of optimal control problems for which all control variables appear linearly. The shooting system has, in the general case, more equations than unknowns and the Gauss-Newton method is used to compute a zero of the shooting function.
This shooting algorithm is locally quadratically convergent, if the derivative of the shooting function is one-to-one at the solution. The main result of this paper is  to show that the latter holds whenever a sufficient condition for weak optimality is satisfied.
We note that this condition is very close to a second order necessary condition.
For the case when the shooting system can be reduced to one having the same number of unknowns and equations (square system), we prove that the mentioned sufficient condition guarantees the stability of the optimal solution under small perturbations and the invertibility of the Jacobian matrix of the shooting function associated to the perturbed problem.
We present numerical tests that validate our method.

\end{abstract}

\if{



 
\textbf{Titre: Un algorithme de tir bien pos\'e pour les probl\`emes de commande optimale avec des arcs singuliers.}

Dans ce travail on pr\'esente pour la premi\`ere fois une condition suffisante pour que l'algorithme de tir soit bien pos\'e quand il est appliqu\'e aux probl\`emes de commande optimale affines dans les commandes.
On commence par \'etudier le cas avec des contraintes initiales-finales sur l'\'etat et commande libre, et en suite on ajoute des contraintes sur la commande.
L'algorithme de tir est bien pos\'e si la d\'eriv\'ee de la fonction de tir associ\'ee est injective dans la solution optimale.
Le r\'esultat principal de cet article montre une condition suffisante pour cette injectivit\'e, qui est tr\`es proche de la condition n\'ecessaire du second ordre.
On montre que cette condition suffisante assure la stabilit\'e de la solution optimale aux petites perturbations et qu'elle garantit aussi que l'algorithme de tir est bien pos\'e pour le probl\`eme perturb\'e.
On pr\'esente des essais num\'eriques qui valident notre m\'ethode.

\textbf{Mots-cl\'es:} Commande optimale, Principe de Pontryaguine, commande singuli\`ere, contraintes sur la commande, algorithme de tir, conditions d'optimalit\'e du second ordre.


\pagestyle{myheadings}
\thispagestyle{plain}
\markboth{M.S. ARONNA, J.F. BONNANS AND P. MARTINON}{A SHOOTING ALGORITHM FOR SINGULAR OPTIMAL CONTROL PROBLEMS}

}\fi

\section{Introduction}

The classical shooting method is used to solve boundary value problems. Hence, it is used to compute the solution of optimal control problems by solving the boundary value problem derived from the Pontryagin Maximum Principle.

Some references can be mentioned regarding the shooting method. The first two works we can find in the literature, dating from years 1956 and 1962, respectively, are Goodman-Lance \cite{GooLan56}  and Morrison et al. \cite{MorRilZan62}. Both present the same method for solving two-point boundary value problems in a general setting, not necessarily related to an optimal control problem. The latter article applies to more general formulations.
The method was studied in detail in Keller's book \cite{Kel68}, and later on Bulirsch \cite{Bul71} applied it to the resolution of optimal control problems.


The case we deal with in this paper, where the shooting method is used to solve optimal control problems with control-affine systems, is treated in, e.g., Maurer \cite{Mau76}, Oberle \cite{OberleThesis,Obe79}, Fraser-Andrews \cite{Fra89b}, Martinon \cite{MartinonThesis} and Vossen \cite{Vos09}. These works provide a series of algorithms and numerical examples with different control structures, but no theoretical foundation is supplied. 
In particular, Vossen \cite{Vos09} deal with a problem in which the control can be written as a function of the state variable, i.e. the control has a feedback representation. He propose an algorithm that involves a finite dimensional optimization problem induced by the switching times. 
Actually, his formulation and the transformation we use for control constrained problems (in Section 8) have similar features, in the sense that both approaches treat the problem by splitting the time interval whenever a switching occurs.
The main difference between Vossen's work and the study here presented is that we treat the general problem (no feedback law is necessary). Furthermore, we justify the well-posedness and the convergence of our algorithm via second order sufficient conditions of the original control problem. 
In some of the just mentioned papers, the control variable had only some of its components entering linearly. This particular structure is studied in more detailed in Aronna \cite{Aro11}, and in the present article we study problems having all affine inputs.



In \cite{BonKup93}, Bonnard and Kupka study the optimal time problem of a generic single-input affine system without control constraints, with fixed initial point and terminal point constrained to a given manifold. 
For this class of problems they establish a link between the injectivity of the shooting function and the optimality of the trajectory by means of the conjugate and focal points theory.
Bonnard et al. \cite{BCT07} provides a survey on
a series of algorithms for the numerical computation of these points, which can be employed to test the injectivity of the shooting function in some cases.
The reader is referred to \cite{BCT07}, Bonnard-Chyba \cite{BonChy} and references therein for further information about this topic.

In addition, Malanowski-Maurer \cite{MalMau96} and Bonnans-Hermant \cite{BonHerIHP09} deal with a problem having mixed control-state and pure state running constraints and satisfying the strong Legendre-Clebsch condition (which does not hold in our affine-input case).
They all establish a link between the invertibility of the Jacobian of the shooting function and some second order sufficient condition for optimality. They provide stability analysis as well.


We start this article by presenting an optimal control problem 
affine in the control, with terminal constraints and free control variables. For this kind of problem, 
we state a set of optimality conditions which is equivalent to the  Pontryagin Maximum Principle.
Afterwards, the second order strengthened generalized Legendre-Clebsch condition is used to eliminate the control variable from the stationarity condition. The resulting set of conditions turns out to be a two-point boundary value problem, i.e. a system of ordinary differential equations having boundary conditions both in the initial and final times.  
We define the shooting function as the mapping that assigns to each estimate of the initial values, the value of the final condition of the corresponding solution. The shooting algorithm consists of approximating a zero of this function. 
In other words, the method finds suitable initial values for which the corresponding solution of the differential equation system satisfies the final conditions.

Since the number of equations happens to be, in general, greater than the number of unknowns, the Gauss-Newton method is a suitable approach for solving this overdetermined system of equations. 
The reader is referred to Dennis \cite{Den77}, Fletcher \cite{Fle80} and Dennis et al. \cite{DenGayWel81} for details and implementations of Gauss-Newton technique.
This method is applicable when the derivative of the shooting function is one-to-one at the solution, and in this case it converges locally quadratically.

The main result of this paper is to provide a sufficient condition for the injectivity of this derivative, and to note that this condition is quite weak since, for qualified problems, it characterizes quadratic growth in the weak sense (see Dmitruk \cite{Dmi77,Dmi87}).
Once the unconstrained case is investigated, we pass to a problem having bounded controls. To treat this case, we perform a transformation yielding a new problem without bounds, we prove that an optimal solution of the original problem is also optimal for the transformed one and we apply our above-mentioned result to this modified formulation.

It is interesting to mention that, by means of the latter result, we can justify, in particular, the invertibility of the Jacobian of the shooting function proposed by Maurer \cite{Mau76}. In this work, Maurer suggested  a method to treat problems having scalar bang-singular-bang solutions and provided a square system of equations (i.e. a system having as many equations as unknowns), meant to be solved by Newton's algorithm. However, the systems that can be encountered in practice may not be square and hence our approach is suitable.


We provide a deeper analysis in the case when the shooting system can be reduced to one having equal number of equations and unknowns. In this framework, we investigate the stability of the optimal solution.
It is shown that the above-mentioned sufficient condition guarantees the stability of the optimal solution under small perturbation of the data, and the invertibility of the Jacobian of the shooting function associated to the perturbed problem.
 Felgenhauer in \cite{Fel11a,Fel11b} provided sufficient conditions for the stability of the structure of the optimal control, but assuming that the perturbed problem had an optimal solution.

Our article is organized as follows. In Section \ref{chap2SecPb}, we present the optimal control problem without bound constraints, for which we provide an optimality system in Section \ref{chap2SecOS}. We give a description of the shooting method in Section \ref{chap2SecShoot}. In Section \ref{chap2SecSOC}, we present a set of second order necessary and sufficient conditions, and the statement of the main result. We introduce a linear quadratic optimal control problem in Section \ref{chap2SecLQ}. 
In Section \ref{chap2SecTransf}, we present a variable transformation relating the shooting system and the optimality system of the linear quadratic problem mentioned above. 
In Section \ref{chap2SecCons}, we deal with the control constrained case.  
A stability analysis for both unconstrained and constrained control cases is provided in Section \ref{chap2SecSta}.
Finally, we present some numerical tests in Section \ref{chap2SecNum}, and we devote Section \ref{chap2SecConclusion} to the conclusions of the article.


\section{Statement of the Problem}\label{chap2SecPb}

Consider the spaces
$\Uspace:=L^{\infty}(0,T;\mathbb{R}^m)$ and
$\mathcal{X}:=W^{1,\infty}(0,T;\mathbb{R}^n),$ as control
and state spaces, respectively. Denote by $u$ and $x$ their elements, respectively.  When needed, put $w=(x,u)$ for a point in the product space $\mathcal{W}:=\mathcal{X}\times \Uspace.$
In this paper, we investigate the optimal control problem
\begin{eqnarray}
 &&\label{chap2cost} J:=\varphi_0(x_0,x_T)\rightarrow \min,\\
 &&\label{chap2stateeq}\dot{x}_t=\sum_{i=0}^mu_{i,t}f_i(x_t),\quad  \ {\rm a.e.}\ {\rm on}\ [0,T],\\
 &&\label{chap2finalcons}  \eta_j(x_0,x_T)=0,\quad \mathrm{for}\ j=1,\hdots,d_{\eta},
\end{eqnarray}
where final time $T$ is fixed, $u_0\equiv 1,$ $f_i:\cR^n\rightarrow \cR^n$ for $i=0,\hdots,m$ and $\eta_j:\cR^{2n}\rightarrow \cR$ for $j=1,\hdots,d_{\eta}.$ 
Assume that data functions $\varphi_0,$ $f_i$ and $\eta_j$ have Lipschitz-continuous second derivatives.
Denote by (P) the problem defined by \eqref{chap2cost}-\eqref{chap2finalcons}.
An element $w\in \mathcal{W}$ satisfying \eqref{chap2stateeq}-\eqref{chap2finalcons} is called a \textit{feasible trajectory}.

Set $\X_*:= W^{1,\infty}(0,T;\cR^{n,*})$ the space of Lipschitz-continuous  functions with values in the $n-$dimensional space $\cR^{n,*}$ of row vectors with real components.
Consider an element $\lambda:=(\beta,p)\in \cR^{d_{\eta},*}\times \X_*$ and define
the \textit{pre-Hamiltonian} function
\benl
H\lam(x,u,t):=p_t\sum_{i=0}^mu_if_i(x),
\eenl
the \textit{initial-final Lagrangian} function
\benl
\ell[\lambda](\zeta_0,\zeta_T)
:=\varphi_0(\zeta_0,\zeta_T)+\sum_{j=1}^{d_{\eta}}\beta_j \eta_j(\zeta_0,\zeta_T),
\eenl
and the \textit{Lagrangian} function
\benl
\L\lam (w):=\ell\lam(x_0,x_T)
+\intT p_t\Big(\sum_{i=0}^m u_{i,t}f_i(x_t)-\dot x_t\Big)\dtt.
\eenl

We study a nominal feasible trajectory $\wh=(\xh,\uh).$
Next, we present a \textit{qualification hypothesis} that is assumed throughout the article. Consider the mapping
\benl 
\ba{rcl}
G\colon  \cR^n\times \U &\rightarrow &\cR^{d_{\eta}}\\
 (x_0,u) & \mapsto &\eta(x_0,x_T),
\ea
\eenl
where $x_T$ is the solution of \eqref{chap2stateeq} associated to $(x_0,u).$

\begin{assumption}
\label{chap2lambdaunique}
The derivative of $G$ at $(\xh_0,\uh)$ is onto.
\end{assumption}

Assumption \ref{chap2lambdaunique} is usually known as \textit{qualification of equality constraints.} 

\begin{definition}
It is said that the trajectory $\wh$ is a {\em weak minimum} of problem (P) iff there exists $\varepsilon>0$ such that $\wh$ is a minimum in the set of feasible trajectories $w=(x,u)\in\W$ satisfying 
\benl
\|x-\xh\|_{\infty}<\varepsilon,\quad  \|u-\uh\|_{\infty}<\varepsilon.
\eenl
\end{definition}

The following first order necessary condition holds for $\wh.$
See the book by Pontryagin et al. \cite{PBGM} for a proof.
\begin{theorem}
\label{chap2Mult}
 If $\wh$ is a weak solution, then there exists an element $\lambda=(\beta,p),$ with $\beta \in \cR^{d_{\eta},*}$ and $p \in \X_*,$ such that $p$ is solution of the {\em costate equation} 
\be\label{chap2costateeq}
-\dot{p}_t =D_xH[\lambda](\xh_t,\uh_t,t),\quad {\rm  a.e.}\ {\rm on}\ [0,T],
\ee
with {\em transversality conditions} 
\begin{align}
 \label{chap2p0} p_0 & =-D_{x_0}\ell[\lambda](\xh_0,\xh_T),\\
 \label{chap2pT} p_T & =D_{x_T}\ell[\lambda](\xh_0,\xh_T),
\end{align}
and the {\em stationarity condition}
\be
\label{chap2stat}
D_uH\lam (\xh_t,\uh_t,t)=0, \quad {\rm  a.e.}\ {\rm  on}\ [0,T],
\ee
is verified.
\end{theorem}

It follows easily that since the pre-Hamiltonian $H$ is affine in all the control variables, \eqref{chap2stat} is equivalent to the  {\em minimum condition}  
\be
\label{chap2maxcond}
H[\lambda](\xh_t,\uh_t,t)=\min_{v\in\cR^m} H[\lambda](\xh_t,v,t), \quad {\rm  a.e.}\ {\rm  on}\ [0,T].
\ee
In order words, the element $(\wh,\lambda)$ in Theorem \ref{chap2Mult} satisfies the \textit{qualified Pontryagin Maximum Principle} and $\lambda$ is a \textit{Pontryagin multiplier.}
On the other hand, it is known that the Assumption \ref{chap2lambdaunique} implies also uniqueness of multiplier. We denote this unique multiplier by $\hat\lambda=(\hat\beta,\ph).$

Let the \textit{switching function} $\Phi:[0,T]\rightarrow \cR^{m,*}$ be defined by
\be 
\label{chap2Phi}
\Phi_t:=D_uH\lamh(\xh_t,\uh_t,t)=(\ph_tf_i(\xh_t))_{i=1}^m.
\ee
Observe that the stationarity condition \eqref{chap2stat} can be written as
\be
\label{chap2Phi0}
\Phi_t=0,\quad {\rm  a.e.}\ {\rm  on}\ [0,T].
\ee


\section{Optimality System}\label{chap2SecOS}

In this section, we present an optimality system, i.e. a set of equations that are necessary for optimality. We obtain this system from the conditions in Theorem \ref{chap2Mult} above and assuming that the \textit{strengthened generalized Legendre-Clebsch condition} (to be defined below) holds.

Observe that, since $H$ is affine in the control, the switching function $\Phi$ introduced in \eqref{chap2Phi} does not depend explicitly on $u.$ 
Let an index $i=1,\hdots,m,$ and $({{\rm  d}^{M_i}}\Phi/{{\rm  d}t^{M_i}})$ be the lowest order derivative of $\Phi$ in which $u_i$ appears with a coefficient that is not identically zero on $]0,T[.$ 
These derivatives of the switching function were used to state necessary condition in Kelley \cite{Kel64}, Goh \cite{Goh66a,Goh66,GohThesis}, Kelley et al. \cite{KelKopMoy67} and Robbins \cite{Rob67}. 
Under the hypothesis that the extremal is normal (as it is the case here by Assumption \ref{chap2lambdaunique}), they proved that the order $M_i$ is even. Hence, the first order derivative $\dot\Phi$ does not depend explicitly on $u$ and a further derivation in time yields 
\be
\label{chap2ddotPhi}
\ddot\Phi_t = 0, \quad {\rm  a.e.}\ {\rm  on}\ [0,T].
\ee
Observe that the latter expression can give explicit information of the control. Actually, in \cite{Goh66a,Goh66,GohThesis,KelKopMoy67,Rob67} it is showed that a necessary condition for weak optimality is that the coefficient of $u$ in \eqref{chap2ddotPhi} satisfies
\be 
\label{chap2genLC}
-\frac{\partial}{\partial u}\ddot\Phi_t \succeq 0,\quad {\rm on}\ [0,T].
\ee
Here, by $X \succeq 0$ we mean that the matrix $X$ is positive semidefinite. The equation \eqref{chap2genLC} is known as {\it generalized Legendre-Clebsch condition.} In order to be able to
express $\uh$ in terms of $(\ph,\xh)$ from \eqref{chap2ddotPhi}, we assume that \eqref{chap2genLC} holds with strict inequality, i.e. we make the following hypothesis.
\begin{assumption}
\label{chap2order}
The {\em strengthened generalized Legendre-Clebsch condition} holds, i.e.
\be 
\label{chap2sgenLC}
-\frac{\partial}{\partial u}\ddot\Phi_t \succ 0,\quad {\rm on}\ [0,T].
\ee
\end{assumption}
Note that function $\ddot\Phi$ is affine in $u,$ and thus $\uh$ can be written in terms of $(\ph,\xh)$ from \eqref{chap2ddotPhi} by inverting the matrix in \eqref{chap2sgenLC}. Due to the regularity hypothesis imposed on the data functions, $\uh$ turns out to be a continuous function of time.
Hence, condition \eqref{chap2ddotPhi} follows from the optimality system and we can use it to compute $\uh$ in view of Assumption \ref{chap2order}. 
In order to guarantee the stationarity condition \eqref{chap2Phi0} we consider the endpoint conditions
\be
\label{chap2endpointPhi}
\Phi_T = 0 ,\   \dot\Phi_0 = 0.
\ee
\begin{remarks}
\label{chap2RemBound}
We could choose another pair of endpoint conditions among the four possible ones: $\Phi_0=0,$ $\Phi_T=0,$ $\dot\Phi_0=0$ and $\dot\Phi_T=0,$ always including at least one of order zero. The choice we made in \eqref{chap2endpointPhi} will simplify the presentation of the results afterwards.
\end{remarks}

\textbf{Notation:} Denote by (OS) the set of equations composed by \eqref{chap2stateeq}-\eqref{chap2finalcons}, \eqref{chap2costateeq}-\eqref{chap2pT}, \eqref{chap2ddotPhi},  \eqref{chap2endpointPhi}, i.e. the system
\be\tag{OS}
\left\{
\begin{split}
\dot{x}_t &= \sum_{i=0}^mu_{i,t}f_i(x_t),\quad  \mathrm{a.e.}\  {\rm on}\ [0,T], \\
\eta_j&(x_0,x_T)=0,\quad \mathrm{for}\ j=1,\hdots,d_{\eta}, \\
 -\dot{p}_t &= D_xH[\lambda](\xh_t,\uh_t,t),\quad \mathrm{a.e.}\  {\rm on}\ [0,T],\\
 p_0 &= -D_{x_0}\ell[\lambda](\xh_0,\xh_T), \quad 
 p_T = D_{x_T}\ell[\lambda](\xh_0,\xh_T),\\
 \ddot\Phi_t &= 0, \quad \mathrm{a.e.}\  {\rm on}\ [0,T],\\
 \Phi_T &= 0 ,\   \dot\Phi_0 = 0.
\end{split}
\right.
\ee

Let us give explicit expressions for $\Phi$ and $\dot\Phi.$ 
Define  
$A \colon \cR^{n+m}\rightarrow \M_{n\times n}(\cR)$ and $B \colon \cR^n\rightarrow \M_{n\times m}(\cR)$ by
\be\label{chap2AB}
A(x,u):=\sum_{i=0}^m u_if_i'(x),\ \ 
B(x)v:=\sum_{i=1}^m v_i f_i(x),
\ee
for every $v\in\cR^m.$ Note that the $i$th. column of $B(x)$ is $f_i(x).$ 
For $(x,u)\in \W$ satisfying \eqref{chap2stateeq}, let $B_1(x_t,u_t) \in \M_{n\times m}(\cR)$ given by 
\be 
\label{chap2B1} 
B_1(x_t,u_t):=A(x_t,u_t)B(x_t)-\ddt {B}(x_t).
\ee 
In view of \eqref{chap2AB} and \eqref{chap2B1}, the expressions in \eqref{chap2endpointPhi} can be rewritten as
\be
\label{chap2formulaePhi}
\Phi_t = p_t B(x_t),\quad 
\dot\Phi_t = -p_t B_1(x_t,u_t).
\ee


\section{Shooting Algorithm}\label{chap2SecShoot}

The aim of this section is to present an appropriated numerical scheme to solve the system (OS). 
For this purpose, define the \textit{shooting function} 
\be
\label{chap2shoot1}
\begin{array}{rl}
\vspace{5pt}\shoot \colon {\rm  D}(\mathcal{S}):=\cR^n\times\cR^{n+d_{\eta},*}&\rightarrow \cR^{d_{\eta}}\times \cR^{2n+2m,*},\\ 
\begin{pmatrix} 
x_0,p_0,\beta 
\end{pmatrix}
 =:\nu
&\mapsto
\mathcal{S}(\nu):=
\begin{pmatrix}
\eta(x_0,x_T)\\
p_0 +D_{x_0}\ell\lam(x_0,x_T)\\
p_T -D_{x_T}\ell\lam(x_0,x_T)\\
p_TB(x_T)\\
p_0 B_1(x_0,u_0)
\end{pmatrix},
\end{array}
\ee
where $(x,u,p)$ is a solution of
 \eqref{chap2stateeq},\eqref{chap2costateeq},\eqref{chap2ddotPhi} corresponding to the initial conditions $(x_0,p_0),$ and with $\lambda:=(\beta,p).$ 
Here, we denote either by $(a_1,a_2)$ or $\begin{pmatrix}a_1\\a_2 \end{pmatrix}$ an element of the product space $A_1\times A_2.$
Note that the control $u$ retrieved from \eqref{chap2ddotPhi} is continuous in time, as we have already pointed out after Assumption \ref{chap2order}. Hence, we can refer to the value $u_0, $ as it is done in the right hand-side of \eqref{chap2shoot1}.
Observe that in a simpler framework having fixed initial state and no final constraints, the shooting function would depend only on $p_0.$ In our case, since the initial state is not fixed and a multiplier associated with the initial-final constraints must be considered,  $\mathcal{S}$ has more independent variables.
Note that solving (OS) consists of finding $\nu \in {\rm  D}(\mathcal{S})$ such that
\be
\label{chap2S=0}
\mathcal{S}(\nu)=0.
\ee
Since the number of equations in \eqref{chap2S=0} is greater than the number of unknowns, the Gauss-Newton method is a suitable approach to solve it. 
This algorithm will solve the equivalent least squares problem
\benl 
\min_{\nu \in {\rm  D}(\mathcal{S})} \left|\mathcal{S}\begin{pmatrix} \nu \end{pmatrix} \right|^2.
\eenl
At each iteration $k,$ given the approximate values $\nu^k,$ it looks for $\Delta^k$ that gives the minimum of the linear approximation of problem
\be\label{chap2minnorm}
\min_{\Delta\in {\rm  D}(\mathcal{S})} \left|\mathcal{S}(\nu^k)+
\mathcal{S}'(\nu^k)\Delta\right|^2.
\ee
Afterwards it updates
\be\label{chap2approxseq}
\nu^{k+1}\leftarrow \nu^k+\Delta^k.
\ee
In order to solve the linear approximation of problem \eqref{chap2minnorm} at each iteration $k,$ we look for $\Delta^k$ in the kernel of the derivative of the objective function, i.e. $\Delta^k$ satisfying
\be
\label{chap2eqS}
\mathcal{S}'(\nu^k)\tras \mathcal{S}'(\nu^k)\Delta^k + \mathcal{S}'(\nu^k)\tras \mathcal{S}(\nu^k)=0.
\ee
Hence, to compute direction $\Delta^k$ the matrix $\mathcal{S}'(\nu^k)\tras \mathcal{S}'(\nu^k)$ must be nonsingular. Thus, {Gauss-Newton method will be applicable, provided that $\mathcal{S}'(\hat\nu)\tras \mathcal{S}'(\hat\nu)$ is invertible,} where $\hat\nu:=(\xh_0,\ph_0,\hat\beta).$ Easily follows that  $\mathcal{S}'(\hat\nu)\tras \mathcal{S}'(\hat\nu)$ is nonsingular if and only if $\mathcal{S}'(\hat\nu)$ is one-to-one.
Summarizing, the \textit{shooting algorithm} we propose here consists of solving the equation \eqref{chap2S=0} by the Gauss-Newton method defined by \eqref{chap2approxseq}-\eqref{chap2eqS}.
 
Since the right hand-side of system \eqref{chap2S=0} is zero, the Gauss-Newton method converges locally quadratically if the function $\shoot$ has Lipschitz-continuous derivative.
The latter holds here given the regularity assumptions on the data functions.
This convergence result is stated in the proposition below. See, e.g., Fletcher \cite{Fle80} or Bonnans \cite{BonOC} for a proof.
\begin{proposition}
\label{chap2Conv}
  If $\mathcal{S}'(\hat\nu)$ is one-to-one then the shooting algorithm is locally quadratically convergent.
\end{proposition}

The main result of this article is to present a condition that guarantees the quadratic convergence of the shooting method near the optimal (local) extremal 
$(\wh,\hat\lambda).$
This condition involves the second variation studied in Dmitruk \cite{Dmi77,Dmi87}, more precisely, 
the sufficient optimality conditions therein presented.

\subsection{Linearization of a Differential Algebraic System}
For the aim of finding an expression of $\mathcal{S}'(\hat\nu),$
we make use of the linearization of (OS) and thus we introduce the following concept.
\begin{definition}[Linearization of a Differential Algebraic System]
Consider a system of differential algebraic equations (DAE) with endpoint conditions
\benl
\dot\zeta_t=\F(\zeta_t,\alpha_t),\quad
0 =\G(\zeta_t,\alpha_t),\quad
0=\I(\zeta_0,\zeta_T),
\eenl
where $\F:\cR^{m+n}\rightarrow \cR^n,$ $\G:\cR^{m+n}\rightarrow \cR^{d_\G}$ and $\I:\cR^{2n}\rightarrow \cR^{d_\I}$ 
are $\C^1-$functions. Let  $(\zeta^0,\alpha^0)$ be a $\C^1-$solution. We call \textit{linearized system} at point $(\zeta^0,\alpha^0)$ the following DAE in the variables $\bar\zeta$ and $\bar\alpha,$
\benl
\dot{\bar\zeta}_t={\rm  Lin}\,\F\mid_{(\zeta^0_t,\alpha^0_t)}(\bar{\zeta}_t,\bar{\alpha}_t),\ 
0={\rm  Lin}\,\G\mid_{(\zeta^0_t,\alpha^0_t)}(\bar{\zeta}_t,\bar{\alpha}_t),\ 
0={\rm  Lin}\,\I\mid_{(\zeta^0_0,\zeta_T^0)}(\bar{\zeta}_0,\bar{\zeta}_T),
\eenl
where
\benl
{\rm  Lin}\,\F\mid_{(\zeta^0_t,\alpha^0_t)}(\bar{\zeta}_t,\bar{\alpha}_t):=\F'(\zeta^0_t,\alpha^0_t)(\bar{\zeta}_t,\bar{\alpha}_t),
\eenl
and the analogous definitions hold for ${\rm  Lin}\,\G$ and ${\rm  Lin}\,\H.$
\if{
\begin{align}
\label{chap2LinF}\dot{\bar\zeta}_t&=\frac{\partial \F}{\partial\zeta}(\zeta^0_t,\alpha^0_t)\bar{\zeta}_t
+\frac{\partial \F}{\partial\alpha}(\zeta^0_t,\alpha^0_t)\bar{\alpha}_t,\\
\label{chap2LinG}0&=\frac{\partial \G}{\partial\zeta}(\zeta^0_t,\alpha^0_t)\bar{\zeta}_t
+\frac{\partial \G}{\partial\alpha}(\zeta^0_t,\alpha^0_t)\bar{\alpha}_t,\\
\label{chap2LinH}0&=\frac{\partial \H}{\partial\zeta_0}(\zeta^0_0,\zeta^0_t)\bar{\zeta}_0
+\frac{\partial \H}{\partial\zeta_T}(\zeta^0_0,\zeta^0_t)\bar{\zeta}_t,\
\end{align}
}\fi
\end{definition}

The technical result below will simplify the computation of the linearization of (OS). Its proof is immediate.

\begin{lemma}
[Commutation of linearization and differentiation]
\label{chap2lindiff}
Given $\G$ and $\F$ as in the previous definition, it holds
\be 
\label{chap2lemmaLineq}
\ddt\,{\rm  Lin}\, \G={\rm  Lin}\,\ddt \G,\quad 
\ddt\,{\rm  Lin}\, \F={\rm  Lin}\,\ddt \F.
\ee 
\end{lemma}

\if{
\proof
Let us compute $\ddt\,{\rm  Lin}\, \G:$
\be\label{chap2ddtLin}
\ba{rl}
&\ddt\,\left[{\rm  Lin}\, \G\mid_{(\zeta^0_t,\alpha^0_t)}(\bar{\zeta}_t,\bar{\alpha}_t)\right]\\
&=[\G_{\zeta\zeta}(\zeta^0_t,\alpha^0_t)\dot\zeta^0_t+
\G_{\zeta\alpha}(\zeta^0_t,\alpha^0_t)\dot\alpha^0_t]\bar\zeta_t
+\G_{\zeta}(\zeta^0_t,\alpha^0_t)\dot{\bar\zeta}_t\\
&+[\G_{\alpha\zeta}(\zeta^0_t,\alpha^0_t)\dot\zeta^0_t+
\G_{\alpha\alpha}(\zeta^0_t,\alpha^0_t)\dot\alpha^0_t]\bar\alpha_t
+\G_{\alpha}(\zeta^0_t,\alpha^0_t)\dot{\bar\alpha}_t.
\ea
\ee
For ${\rm  Lin}\,\ddt \G$ we get
\be\label{chap2Linddt}
\ba{rl}
&{\rm  Lin}\mid_{(\zeta^0_t,\alpha^0_t)}\,\ddt \G(\zeta_t,\alpha_t)\\
&={\rm  Lin}\mid_{(\zeta^0_t,\alpha^0_t)}
[\G_{\zeta}(\zeta_t,\alpha_t)\dot\zeta_t+
\G_{\alpha}(\zeta_t,\alpha_t)\dot\alpha_t]\\
&=[\G_{\zeta\zeta}(\zeta^0_t,\alpha^0_t)\bar\zeta_t
+\G_{\zeta\alpha}(\zeta^0_t,\alpha^0_t)\bar\alpha_t]\dot{\zeta}^0_t
+\G_{\zeta}(\zeta^0_t,\alpha^0_t)\dot{\bar\zeta}_t\\
&+[\G_{\alpha\zeta}(\zeta^0_t,\alpha^0_t)\bar\zeta_t
+\G_{\alpha\alpha}(\zeta^0_t,\alpha^0_t)\bar\alpha_t]
\dot\alpha^0_t
+\G_{\alpha}(\zeta^0_t,\alpha^0_t)\dot{\bar\alpha}_t.
\ea
\ee
By \eqref{chap2ddtLin} and \eqref{chap2Linddt} we get \eqref{chap2lemmaLineq}.
\qed
}\fi


\subsection{Linearized Optimality System}\label{chap2SecLS}

In the sequel, whenever the argument of functions $A,B,B_1,$ etc. is omitted, assume that they are evaluated at the reference extremal $(\wh,\hat\lambda).$ 
Define the $m\times n-$matrix $C,$ the $n\times n-$matrix $Q$ and the $m\times n-$matrix $M$ by
\be\label{chap2CQ}
C:=H_{ux},
\quad 
Q:=H_{xx},
\quad 
M:=B^\top Q- \dot C-C A.
\ee
Note that the $i$th. row of matrix $C$ is the function $pf_i',$ for $i=1,\hdots,m.$
Denote with $(z,v,\bar\lambda:=(\bar\beta,q))$ the linearized variable $(x,u,\lambda=(\beta,p)).$
In view of equations \eqref{chap2formulaePhi} and \eqref{chap2CQ} we can write
\be
\label{chap2LinPhi}
{\rm  Lin}\ \Phi_t=q_tB_t + z\tras_tC\tras_t.
\ee
The linearization of system (OS) at point $(\xh,\uh,\hat\lambda)$ consists of the \textit{linearized state equation}
\be
\label{chap2eqz}
 \dot z_t = A_tz_t+B_tv_t,\quad {\rm a.e.}\ {\rm on}\ [0,T],
\ee
with endpoint conditions
\be 
\label{chap2Linfinal} 
 0 = D\eta(\xh_0,\xh_T)(z_0,z_T),
\ee
the linearized costate equation
\be
\label{chap2eqq}
 -\dot q_t = q_tA_t + z_t\tras Q_t+v_t\tras C_t,\quad {\rm  a.e.}\ {\rm  on}\ [0,T],
\ee
with endpoint conditions
\begin{align}
\label{chap2condq0} 
  q_0 &= -\left[z_0\tras D^2_{x_0^2}\ell+z_T\tras D^2_{x_0x_T}\ell 
+\sum_{j=1}^{d_{\eta}}{\bar\beta}_jD_{x_0}\eta_j\right]_{(\xh_0,\xh_T)},\\
\label{chap2condqT}
q_T&=\left[z_T\tras D^2_{x_T^2}\ell+z_0\tras D^2_{x_0x_T}\ell 
+\sum_{j=1}^{d_{\eta}}{\bar\beta}_jD_{x_T}\eta_j\right]_{(\xh_0,\xh_T)},
\end{align}
and the algebraic equations
\begin{align}
 \label{chap2LinddotHu}  
0&={\rm  Lin}\ \ddot\Phi=-\frac{{\rm d}^2}{{\rm d}t^2}(qB+Cz),\quad {\rm  a.e.}\ {\rm  on}\ [0,T],\\ 
\label{chap2LinPhiT}
0&= {\rm  Lin}\ \Phi_T=q_TB_T+C_T z_T,\\
\label{chap2LindotPhi0}
0&= {\rm  Lin}\ \dot\Phi_0=-\ddt(qB+Cz)_{t=0}.
\end{align}
Here we used equation \eqref{chap2LinPhi} and commutation property of Lemma \ref{chap2lindiff} to write \eqref{chap2LinddotHu} and \eqref{chap2LindotPhi0}. 
Observe that \eqref{chap2LinddotHu}-\eqref{chap2LindotPhi0} and Lemma \ref{chap2lindiff} yield
\be
\label{chap2LinPhi0}
0={\rm  Lin}\ \Phi_t=q_tB_t+z_t\tras C_t\tras,\quad {\rm  on}\ [0,T],
\ee
and
\benl
0={\rm  Lin}\ \dot\Phi_t=-qB_1-z\tras M\tras + v\tras(-CB+B\tras C\tras),
\quad {\rm a.e.}\ {\rm  on}\ [0,T].
\eenl
By means of Theorem \ref{chap2thGoh}, to be stated in Section \ref{chap2SecSOC} afterwards, we can see that the coefficient of $v$ in previous expression vanishes, and hence,
\be
\label{chap2LindotPhi0}
0={\rm  Lin}\ \dot\Phi_t=-qB_1-z\tras M\tras ,
\quad 
{\rm  on}\ [0,T].
\ee
Note that both equations \eqref{chap2LinPhi0} and \eqref{chap2LindotPhi0} hold everywhere on $[0,T]$ since all the involved functions are continuous in time.

\noindent\textbf{Notation:} denote by (LS) the set of equations \eqref{chap2eqz}-\eqref{chap2LindotPhi0}.

Once we have computed the linearized system (LS), we can write the derivative of $\shoot$ in the direction  $\bar\nu:=\begin{pmatrix} z_0, q_0, \bar\beta \end{pmatrix}$ as follows.
\benl
\shoot'
(\hat\nu)\bar\nu
=
\begin{pmatrix}
\vspace{4pt} D\eta(\xh_0,\xh_T)(z_0,z_T)\\
\vspace{4pt}q_0+\left[z_0\tras D^2_{x_0^2}\ell+z_T\tras D^2_{x_0x_T}\ell 
+\sum_{j=1}^{d_{\eta}}\bar{\beta}_jD_{x_0}\eta_j\right]_{(\xh_0,\xh_T)}\\
\vspace{4pt} q_T-\left[z_T\tras D^2_{x_T^2}\ell+z_0\tras D^2_{x_0x_T}\ell 
+\sum_{j=1}^{d_{\eta}}\bar{\beta}_jD_{x_T}\eta_j\right]_{(\xh_0,\xh_T)}\\
\vspace{4pt} q_TB_T+z_T\tras C_T\tras\\
q_0B_{1,0}+z_0\tras M\tras_0
\end{pmatrix},
\eenl
where  $(v,z,q)$ is the solution of \eqref{chap2eqz},\eqref{chap2eqq},\eqref{chap2LinddotHu} associated with the initial condition $(z_0,q_0)$ and the multiplier $\bar\beta.$
Thus, we get the property below.
\begin{proposition}
\label{chap2S'}
$\shoot'(\hat\nu)$  is one-to-one if the only solution of \eqref{chap2eqz}-\eqref{chap2eqq},\eqref{chap2LinddotHu} is $(v,z,q)=0.$ 
\end{proposition}


\section{Second Order Optimality Conditions}\label{chap2SecSOC}

In this section we summarize a set of second order necessary and sufficient conditions. At the end of the section we state a sufficient condition for the local quadratic convergence of the shooting algorithm presented in Section \ref{chap2SecShoot}. The latter is the main result of this article.

Recall the matrices $C$ and $Q$ defined in \eqref{chap2CQ}, and
the space $\W$ given at the beginning of Section \ref{chap2SecPb}.
Consider the quadratic mapping on $\W,$  
\be\label{chap2Omega}
\Omega(z,v):=\half D^2\ell\,(z_0,z_T)^2
+\half\int_0^T \left[z\tras Q z+2 v\tras C z\right]\dtt.
\ee
It is a well-known result that for each $(z,v)\in \W,$
\be 
\half D^2\L\, (z,v)^2 =\Omega(z,v).
\ee
We next recall the classical second order necessary condition for optimality that states that the second variation of the Lagrangian function is nonnegative on the critical cone. In our case, the \textit{critical cone} is given by
\be 
\label{chap2coneC}
\C:=\{(z,v)\in \W:\,\text{\eqref{chap2eqz}-\eqref{chap2Linfinal} hold} \},
\ee
and the second order optimality condition is as follows.
\begin{theorem}[Second order necessary optimality condition]
\label{chap2NC}
If $\wh$ is a weak minimum of (P), then
\be
\label{chap2NCeq}
\Omega(z,v)\geq 0,\quad {\rm for}\ {\rm all}\ (z,v)\in\C.
\ee
\end{theorem}
A proof of previous theorem can be found in, e.g.,  Levitin, Milyutin and Osmolovskii \cite{LMO}.

In the sequel we present a necessary condition due to Goh \cite{Goh66} and a second order sufficient condition from Dmitruk \cite{Dmi77}.
The idea behind these results lies on the following observation. 
Note that the mapping $\Omega$ in the necessary condition of Theorem 3 does not contain a quadratic term on $v$ (since $H_{uu}\equiv 0$). Hence, one cannot deduce a necessary condition from \eqref{chap2NCeq} in terms of the positive semidefiniteness of some matrix, as it is done in the non-affine control case. Furthermore, one cannot attempt to obtain a sufficient condition by strengthening the inequality \eqref{chap2NCeq}. In order to overcome this inconvenience,
Goh introduced a change of variables in \cite{Goh66a} and employed it to derive necessary conditions in \cite{Goh66a,Goh66}. Afterwards, Dmitruk in \cite{Dmi77} stated a second order sufficient condition in terms of the coercivity of $\Omega$ in the corresponding transformed space of variables. 
Let us give the details of this transformation and the transformed second variation.
Given $(z,v)\in \W,$ define
\be
 \label{chap2ydef} y_t:=\int_0^t v_s{\rm  d}s,\quad
  \xi_t:=z_t-B(\xh_t)y_t.
\ee
This change of variables can be performed in any linear system of differential equations, and it is known as \textit{Goh's transformation}.

We aim to perform Goh's transformation in \eqref{chap2Omega}. To this end, 
 consider the spaces $\U_2:=L_2(0,T;\cR^m),$  $\X_2:=W_2^1(0,T;\cR^n),$ 
the function $g\colon \cR^{2n+m} \rightarrow \cR,$ with
\benl
g(\zeta_0,\zeta_T,h):=
D^2\ell\,(\zeta_0,\zeta_T+B_Th)^2
+
h^\top C_T(2\zeta_T+B_Th),
\eenl
and the quadratic mapping ${\bar\Omega}\colon\X_2\times \U_2\times \cR^m \rightarrow \cR,$ given by
\be\label{chap2OmegaP}
(\xi,y,h) \mapsto \bar\Omega(\xi,y,h):=
\half g(\xi_0,\xi_T,h)
+\half \ds\int_0^T\{\xi\tras Q\xi+ 2y\tras M\xi
+ y\tras Ry\}\dtt,
\ee
where the involved matrices where introduced in \eqref{chap2AB}, \eqref{chap2CQ} and \eqref{chap2R}.

The following result is due to Goh \cite{Goh66} and it is a nontrivial consequence of Theorem \ref{chap2NC}. Define the $m\times m-$matrix
\be\label{chap2R}
R:=B^\top Q B-C B_1-(CB_1)\tras -\ddt (C B).
\ee
\begin{theorem}[Goh's Necessary Condition]
\label{chap2thGoh}
If $\wh$ is a weak minimum of (P), then
\be
\label{chap2CBsym} 
C B\ {\rm is}\  {\rm symmetric.}
\ee
Furthermore, 
\benl
\label{chap2Omegas=}
\Omega(z,v)=\bar\Omega(\xi,y,y_T),
\eenl
whenever $(z,v)\in \W$ and $(\xi,y,y_T)\in \X\times \Y\times\cR^m$ satisfy \eqref{chap2ydef}, and the coefficient of the quadratic term in $y$ in $\bar\Omega$ satisfies
\be
\label{chap2Rsucceq0} 
R\succeq 0.
\ee
\end{theorem}

Theorem  \ref{chap2thGoh} was first proved by Goh in \cite{Goh66}, but the notation used here differs a lot from the one employed by Goh.  For a proof of this Theorem in the present notation the reader can consult the recent article by Aronna et el. \cite{ABDL11}.

\begin{remarks}\label{remLC}
Following the analysis in Goh \cite{GohThesis} (see also Kelley et al. \cite{KelKopMoy67}, Bell-Jacobson \cite{BelJac} and references therein), we have that
$R=-\ds\frac{\partial}{\partial u}\ddot\Phi_t,$ and hence, conditions \eqref{chap2Rsucceq0} and \eqref{chap2genLC} are equivalent.
\end{remarks}

\begin{remarks}
Observe that \eqref{chap2CBsym} is equivalent to $pf_i'f_j=pf_j'f_i,$ for every pair $i,j=1,\hdots,m.$ These identities can be written in terms of Lie brackets as
\benl 
p[f_i,f_j]=0,\quad {\rm for}\ i,j=1,\hdots,m.
\eenl
Here $[g,h]$ denotes the \textit{Lie bracket} of two smooth vector fields $g,h\colon\cR^n\rightarrow\cR^n$  and it defined by
\benl 
[g,h](x):=g'(x)h(x)-h'(x)g(x).
\eenl
Note that \eqref{chap2CBsym} implies, in view of \eqref{chap2R}, that $R$ is symmetric.
\if{The components of matrix $R$ can be written as 
\benl 
R_{ij}=p[f_i,[f_j,f_0]],
\eenl
and hence, its symmetry implies
\benl 
p[f_i,[f_j,f_0]]=p[f_j,[f_i,f_0]],\quad {\rm for}\ i,j=1,\hdots,m.
\een}\fi
The latter expressions involving Lie brackets can be often found in the literature.
\end{remarks}

Define the \textit{order function} $\gamma\colon \cR^n \times \U_2\times\cR^m \rightarrow \cR$  as
\benl
\gamma(\xi_0,y,h):=|\xi_0|^2+\intT |y_t|^2\dtt+|h|^2.
\eenl
We call $(\delta x,v)\in \W$ a \textit{feasible variation} for $\wh$ if $(\xh+\delta x,\uh+v)$ satisfies \eqref{chap2stateeq}-\eqref{chap2finalcons}.
\begin{definition}
We say that $\wh$ satisfies the {\em $\gamma-$growth condition in the weak sense} iff there exists $\rho>0$ such that, for every sequence of feasible variations $\{(\delta x^k,v^k)\}$ converging to 0 in $\W,$
\benl
J(\uh+v^k)-J(\uh)\geq \rho\gamma(\xi_0^k,y^k,y^{k}_T),
\eenl
holds for big enough $k,$ where $y^k_t :=\int_0^t v^k_s{\rm  d}s,$
and $\xi^k$ is given by \eqref{chap2ydef}.
\end{definition}

In the previous definition, given that $(\delta x^k,v^k)$ is a feasible variation for each $k,$ the sequence $\{(\delta x^k,v^k)\}$ goes to 0 in $\W$ if and only if $\{ v^k\}$ goes to 0 in $\U.$

Observe that, if $(z,v)\in \W$ satisfies \eqref{chap2eqz}-\eqref{chap2Linfinal},  then $(\xi,y,h:=y_T)$ given by transformation  \eqref{chap2ydef} verifies
\begin{align}
\label{chap2tlineareq} 
&\dot\xi=A\xi+B_1y,\\
 \label{chap2tLinfinal} & D\eta(\xh_0,\xh_T)(\xi_0,\xi_T+B_Th)=0.
\end{align}
Set the \textit{transformed critical cone}
\benl
\P_2:= \left\{
(\xi,y,h)\in \X_2\times\Uspace_2\times \cR^m:
\text{\eqref{chap2tlineareq}-\eqref{chap2tLinfinal}}\ {\rm  hold}\right\}.
\eenl

The following is an immediate consequence of the sufficient condition established in Dmitruk \cite{Dmi77} (or \cite[Theorem 3.1]{Dmi87}).

\begin{theorem}
\label{chap2sc2}
The trajectory $\wh$ is a weak minimum of (P) satisfying $\gamma-$ growth condition in the weak sense if and only if  \eqref{chap2CBsym} holds and there exists $\rho\gr 0$ such that
\be\label{chap2unifpos}
\bar\Omega(\xi,y,h)\geq \rho\gamma(\xi_0,y,h),\quad {\rm  on}\ \P_2.
\ee
\end{theorem}

The result presented in \cite{Dmi77} applies to a more general case having finitely many equalities and  inequalities constraints on the initial and final state, and a set of multipliers consisting possibly of more than one element. 

\begin{remarks} 
\label{chap2RemStrengthLC}
If \eqref{chap2unifpos} holds, then necessarily 
\be
\label{chap2StrengthLC}
R \succeq \rho\, I_m,
\ee
where $I_m$ represents the $m\times m-$identity matrix. Hence, in view of Remark \ref{remLC}, the uniform positivity \eqref{chap2unifpos} implies the strengthened Legendre-Clebsch condition in Assumption \ref{chap2order}.
\end{remarks}

\begin{theorem}\label{chap2wp} 
If $\wh$ is a weak minimum of (P) satisfying 
\eqref{chap2unifpos}, then the shooting algorithm 
is locally quadratically convergent.
\end{theorem}

We present the proof of previous theorem at the end of Section \ref{chap2SecTransf}. 

\begin{remarks}
It is interesting to observe that condition \eqref{chap2unifpos} is a quite weak assumption in the sense that 
it is necessary for $\gamma-$growth and its corresponding relaxed condition \eqref{chap2NCeq} holds necessarily for every weak minimum.
\end{remarks}

\begin{remarks}[Verification of \eqref{chap2unifpos}]
\label{chap2RemVerification}
The sufficient condition in \eqref{chap2unifpos} can be sometimes checked analytically. 
On the other hand, when the initial point $\xi_0$ is fixed, it can be characterized by a Riccati-type equation and/or the nonexistence of a focal point as it was established in Zeidan \cite{Zei92}.
Furthermore, under certain hypotheses, the condition \eqref{chap2unifpos} can be verified numerically as proposed in \cite{COTCOT} by Bonnard, Caillau and Tr\'elat (see also the survey in \cite{BCT07}).
\end{remarks}


\section{Corresponding Linear-quadratic Problem}\label{chap2SecLQ}

In this section, we study the linear-quadratic problem (LQ) given by
\begin{eqnarray}
&&\label{chap2tcost}{\bar\Omega}(\xi,y,h_T)\rightarrow \min,\\
&&\text{\eqref{chap2tlineareq}-\eqref{chap2tLinfinal}},\\
&&\label{chap2heq}\dot{h}=0,\quad h_0\ {\rm  free}.
\end{eqnarray}
Here, $y$ is the control, $\xi$ and $h$ are the state variables.
Note that, if condition \eqref{chap2unifpos} holds, then (LQ) has a unique optimal solution $(\xi,y,h)=0.$
Furthermore, recall that \eqref{chap2unifpos} yields \eqref{chap2StrengthLC} as  it was said in Remark \ref{chap2RemStrengthLC}. In other words, \eqref{chap2unifpos} implies that the strengthened Legendre-Clebsch condition is verified at $(\xi,y,h)=0.$ 
Hence, the unique local optimal solution of (LQ) is characterized by the first optimality system, that we denote afterwards by (LQS).
In Section \ref{chap2SecTransf}, we present a one-to-one linear mapping that transforms each solution of (LS)  (introduced in paragraph \ref{chap2SecLS}) into a solution of this new optimality system (LQS). Theorem \ref{chap2wp} will follow.

Denote by $\chi$ and $\chi_h$  the costate variables corresponding to $\xi$ and $h,$ respectively; and by 
$\beta^{LQ}$ the multiplier associated to the initial-final linearized state constraint \eqref{chap2tLinfinal}.
 Note that the qualification hypothesis in Assumption \ref{chap2lambdaunique}
implies that $\{D\eta_j(\xh_0,\xh_T)\}_{j=1}^{d_{\eta}}$ are linearly independent. 
Hence any weak solution $(\xi,y,h)$ of (LQ) has a unique associated multiplier $\lambda^{LQ}:=(\chi,\chi_h, \beta^{LQ})$ solution of the system that we describe next.
The pre-Hamiltonian of (LQ) is
\benl
\mathcal{H}\lamLQ(\xi,y):=\chi(A\xi+B_1y)+\half (\xi^\top Q \xi +2 y^\top M \xi +  y^\top R y).
\eenl
Observe that $\H$ does not depend on $h$ since the latter has zero dynamics and does not appear in the running cost. 
The endpoint Lagrangian is given by
\benl
\ell^{LQ}\lamLQ(\xi_0,\xi_T,h_T):=
\half g(\xi_0,\xi_T,h_T)
+\sum_{j=1}^{d_{\eta}}\beta_j^{LQ} D\eta_j(\xi_0,\xi_T+B_Th_T).
\eenl
The costate equation for $\chi$ is
\be
\label{chap2chi}
-\dot\chi  =D_{\xi} \mathcal{H}\lamLQ = \chi A + \xi^\top Q + y^\top M,
\ee
with endpoint conditions
\begin{align} 
\label{chap2chi0}
&\ba{rl}
\chi_0=&
-D_{\xi_0}\ell^{LQ}\lamLQ\\
=&
-\left[ \xi_0\tras D_{x_0^2}^2\ell+(\xi_T+B_Th)\tras D_{x_0x_T}^2\ell
+\sum_{j=1}^{d_{\eta}}\beta_j^{LQ} D_{x_0}\eta_j \right],
\ea
\\
\label{chap2chiT}
&\ba{rl}
\chi_T=&
D_{\xi_T}\ell^{LQ}\lamLQ\\
=&
\xi_0\tras D_{x_0x_T}^2\ell
+(\xi_T+B_Th)\tras D^2_{x_T^2}\ell 
+h\tras C_T
+\sum_{j=1}^{d_{\eta}}\beta_j^{LQ} D_{x_T}\eta_j.
\ea
\end{align}
For costate variable $\chi_h$ we get the equation
\be
\label{chap2chihT0}
\dot\chi_{h}=0,\quad
\chi_{h,0}=0,\quad
\chi_{h,T}
=D_h\ell^{LQ}\lamLQ.
\ee
Hence, $\chi_h\equiv 0$ and thus, the last identity in \eqref{chap2chihT0} yields
\be
\label{chap2chihT}
0= \xi_0\tras D^2_{x_0x_T}\ell B_T+(\xi_T+B_Th)\tras(D^2_{x_T^2}\ell B_T+ C_T\tras) 
+\sum_{j=1}^{d_{\eta}}\beta_j^{LQ} D_{x_T}\eta_jB_T.
\ee
The {stationarity with respect to the new control $y$} implies
\be\label{chap2Hy}
0=D_{y} \H=\chi B_1+\xi\tras M\tras+y\tras R.
\ee
\noindent\textbf{Notation:} Denote by (LQS) the set of equations consisting of \eqref{chap2tlineareq}-\eqref{chap2tLinfinal},
\eqref{chap2heq},\eqref{chap2chi}-\eqref{chap2chiT},\eqref{chap2chihT} and \eqref{chap2Hy}, i.e. (LQS) is the system
\benl
\left\{
\begin{split}
 \dot\xi&=A\xi+B_1y,\\
 D&\eta(\xh_0,\xh_T)(\xi_0,\xi_T+B_Th)=0,\\
 \dot{h}&=0,\\
-\dot\chi & =D_{\xi} \mathcal{H}\lamLQ = \chi A + \xi^\top Q + y^\top M,\\
\chi_0&= -\left[ \xi_0\tras D_{x_0^2}^2\ell+(\xi_T+B_Th)\tras D_{x_0x_T}^2\ell
+\sum_{j=1}^{d_{\eta}}\beta_j^{LQ} D_{x_0}\eta_j \right],\\
\chi_T&=\xi_0\tras D_{x_0x_T}^2\ell
+(\xi_T+B_Th)\tras D^2_{x_T^2}\ell 
+h\tras C_T
+\sum_{j=1}^{d_{\eta}}\beta_j^{LQ} D_{x_T}\eta_j,\\
0 &= \xi_0\tras D^2_{x_0x_T}\ell B_T+(\xi_T+B_Th)\tras(D^2_{x_T^2}\ell B_T+ C_T\tras) 
+\sum_{j=1}^{d_{\eta}}\beta_j^{LQ} D_{x_T}\eta_jB_T,\\
0&=\chi B_1+\xi\tras M\tras+y\tras R.
\end{split}
\right.
\eenl
Note that (LQS) is a {first order optimality system} for problem \eqref{chap2tcost}-\eqref{chap2heq}.


\section{The Transformation}\label{chap2SecTransf}

In this section we show how to transform a solution of (LS) into a solution of (LQS) via a one-to-one linear mapping. Given $(z,v,q,\bar\beta)\in\X\times \U\times \X_*\times \cR^{d_{\eta},*},$ define
\be\label{chap2transf}
y_t:=\int_0^t v_sds,\ \xi:=z-By,\ \chi:=q+y\tras C,\ \chi_h:=0,\ h:=y_T,\ 
\beta_j^{LQ}:=\bar\beta_j.
\ee
The next Lemma shows that 
the point $(\xi,y,h,\chi,\chi_h,\beta^{LQ})$ is solution of (LQS) provided that $(z,v,q,\bar\beta)$ is solution of (LS).
\begin{lemma}\label{chap2lemmatransf}
The one-to-one linear mapping defined by \eqref{chap2transf} converts each solution of (LS) into a solution of (LQS).
\end{lemma}
\proof
Let $(z,v,q,\bar\beta)$ be a solution of (LS), and set
$(\xi,y,\chi,\beta^{LQ})$  by \eqref{chap2transf}.

\noindent\textbf{Part I.} We shall prove that $(\xi,y,\chi,\beta^{LQ})$ satisfies conditions \eqref{chap2tlineareq} and \eqref{chap2tLinfinal}.
Equation  \eqref{chap2tlineareq} follows by differentiating expression of $\xi$ in \eqref{chap2transf}, and equation \eqref{chap2tLinfinal} follows from  \eqref{chap2Linfinal}.

\noindent\textbf{Part II.} We shall prove that $(\xi,y,\chi,\beta^{LQ})$ verifies \eqref{chap2chi}-\eqref{chap2chiT} and \eqref{chap2chihT}.
Differentiate $\chi$ in \eqref{chap2transf}, use equations \eqref{chap2eqq} and \eqref{chap2transf}, recall the definition of $M$ in \eqref{chap2CQ} and obtain
\benl
\ba{rl}
-\dot\chi
=& -\dot q-v\tras C-y\tras \dot C = qA + z\tras Q-y\tras \dot C\\
=& \chi A+\xi\tras Q+y\tras (-CA +B\tras Q-\dot C) = \chi A + \xi\tras Q+y\tras M.
\ea
\eenl
Hence \eqref{chap2chi} holds.
Equations \eqref{chap2chi0} and \eqref{chap2chiT} follow from \eqref{chap2condq0} and \eqref{chap2condqT}.
Combine \eqref{chap2condqT} and \eqref{chap2LinPhiT} to get
\benl
\ba{rl}
 0
=&\, q_TB_T + z_T\tras C_T\tras\\
=& \left[z_T\tras D^2_{x_T^2}\ell +z_0\tras D^2_{x_0x_T}\ell  
+\sum_{j=1}^{d_{\eta}}\bar{\beta}_jD_{x_T}\eta_j\right]_{(\xh_0,\xh_T)} B_T+ z_T\tras C_T\tras.
\ea
\eenl
Performing transformation \eqref{chap2transf} in the previous equation  yields \eqref{chap2chihT}.

\noindent\textbf{Part III.} We shall prove that \eqref{chap2Hy} holds. 
Differentiating \eqref{chap2LinPhi0} we get
\benl
0=\ddt {\rm  Lin}\ \Phi=\ddt(qB+z\tras C\tras).
\eenl
Consequently, by \eqref{chap2eqz} and \eqref{chap2eqq},
\be\label{chap2Hy1}
0= -(qA+z\tras Q + v\tras C)B+q\dot B+(z\tras A\tras + v\tras B\tras)C\tras +z\tras \dot C\tras,
\ee
where the coefficient of $v$ vanishes in view of \eqref{chap2CBsym}.
Recall \eqref{chap2B1} and \eqref{chap2CQ}.
Performing transformation \eqref{chap2transf} in \eqref{chap2Hy1} leads to
\benl
0
= -\chi B_1-\xi\tras M\tras+y\tras(CB_1 -B\tras QB+B\tras A\tras C\tras + B\tras\dot C\tras).
\eenl
Equation \eqref{chap2Hy} follows from \eqref{chap2R} and condition \eqref{chap2CBsym}.

Parts I, II and III show that $(\xi,y,\chi,\beta^{LQ})$ is a solution of (LQS), and hence, the result follows.

\qed

\begin{remarks} 
Observe that the unique assumption we needed in previous proof was Goh's condition \eqref{chap2CBsym} that follows from the weak optimality of $\wh.$
\end{remarks}

\proof
\ [of Theorem \ref{chap2wp}]
 We shall prove that \eqref{chap2unifpos} implies that $\shoot'(\hat\nu)$ is one-to-one.
Take $(z,v,q,\bar\beta)$ a solution of (LS), and let $(\xi,y,\chi,\chi_h,\beta^{LQ})$ be defined by \eqref{chap2transf}, that we know by Lemma \ref{chap2lemmatransf} is solution of (LQS). 
As it has been already pointed out at the beginning of Section \ref{chap2SecLQ}, condition \eqref{chap2unifpos} implies that the unique solution of (LQS) is 0.
Hence $(\xi,y,\chi,\chi_h,\beta^{LQ})=0$ and thus $(z,v,q,\bar\beta)=0.$ 
Conclude that the unique solution of (LS) is 0. 
The latter assertion implies, in view of Proposition \ref{chap2S'}, that 
$\shoot'(\hat\nu)$ is one-to-one.
The result follows from Proposition \ref{chap2Conv}.
\qed



\section{Control Constrained Case}\label{chap2SecCons}

In this section, we add the following bounds to the control variables
\be
\label{chap2contcons}
0\leq u_{i,t}\leq 1,\quad {\rm  a.e.}\ {\rm on}\ [0,T],\ {\rm  for}\ i=1,\hdots,m.
\ee
Denote with (CP) the problem given by \eqref{chap2cost}-\eqref{chap2finalcons} and \eqref{chap2contcons}.
\begin{definition}
A feasible trajectory $\wh\in\W$ is  a {\em Pontryagin minimum} of (CP) iff for any positive $N,$ there exists $\varepsilon_N>0$ such that $\wh$ is a minimum in the set of feasible trajectories $w=(x,u)\in\W$ satisfying 
\benl
\|x-\xh\|_{\infty}<\varepsilon_N,\ \|u-\uh\|_{1}<\varepsilon_N,\ \|u-\uh\|_{\infty}<N.
\eenl
\end{definition}
Given $i=1,\hdots,m,$ we say that $\uh_i$ has a \textit{bang arc} on an interval $I\subset [0,T]$ iff $\uh_{i,t}=0$ a.e. on $I,$ or $\uh_{i,t}=1$ a.e. on $I,$ and it has a \textit{singular arc} iff $0\mi  \uh_{i,t} \mi 1$ a.e. on $I.$

\begin{assumption}\label{chap2geohyp}
Each component $\uh_i$ is a finite concatenation of bang and singular arcs.
\end{assumption}

A time $t\in ]0,T[$ is called \textit{switching time} iff there exists an index $1\leq i\leq m$ such that $\uh_i$ switches at time $t$  from singular to bang, or vice versa, or from one bound in \eqref{chap2contcons} to the other.

\begin{remarks}
Assumption \ref{chap2geohyp} rules out the solutions having an infinite number of switchings in a bounded interval. This behavior is usually known as Fuller's phenomenon (see Fuller  \cite{Ful63}).
Many examples can be encountered satisfying Assumption \ref{chap2geohyp} as is the case of the three problems presented in Section \ref{chap2SecNum}.
\end{remarks} 

With the purpose of solving (CP) numerically, we assume that the structure of the concatenation of bang and singular arcs of the optimal solution $\wh$ and an approximation of its switching times are known. 
This initial guess can be obtained, for instance, by solving the nonlinear problem resulting from the discretization of the optimality conditions or by a continuation method. See Betts \cite{Bet98} or Biegler \cite{Bie10} for a detailed survey and description of numerical methods for nonlinear programming problems. For the continuation method the reader is referred to Martinon \cite{MartinonThesis}. 

This section is organized as follows. From (CP) and the known structure of $\uh$ and its switching times we create a new problem that we denote by (TP). Afterwards we prove that we can transform $\wh$ into a weak solution $\Wh$ of (TP). Finally we conclude that if $\Wh$ satisfies the coercivity condition \eqref{chap2unifpos}, then the shooting method for problem (TP)  converges locally quadratically.
In practice, the procedure will be as follows: obtain somehow the structure of the optimal solution of (CP), create problem (TP), solve (TP) numerically obtaining $\Wh,$ and finally transform $\Wh$ to find $\wh.$

Next, we present the transformed problem.

\begin{assumption}\label{chap2disc}
 Assume that each time a control $\uh_i$ switches from bang to singular or vice versa, there is a discontinuity of first kind.
\end{assumption}

Here, by \textit{discontinuity of first kind} we mean that each component of $\uh$ has a finite nonzero jump at the switching times,
and the left and right limits exist.

By Assumption \ref{chap2geohyp} the set of switching times is finite. Consider the partition of $[0,T]$ induced by the switching times:
\benl
\{0=:\Th_0\mi \Th_1\mi \hdots \mi \Th_{N-1}\mi \Th_N:=T\}.
\eenl
Set $\hat{I}_k:=[\Th_{k-1},\Th_k],$ and define for $k=1,\hdots,N,$
\begin{align*}
S_k&:= \{ 1\leq i\leq m:\ \uh_i\ {\rm  is}\ {\rm singular}\ {\rm  on}\ \hat{I}_k\}, \\
E_k&:= \{ 1\leq i\leq m:\ \uh_i=0\ {\rm  a.e.}\ {\rm  on}\ \hat{I}_k \},\\
N_k&:= 
\{ 1\leq i\leq m:\ \uh_i=1\ {\rm  a.e.}\ {\rm  on}\ \hat{I}_k\}.
\end{align*}
Clearly $S_k \cup E_k \cup N_k=\{1,\hdots,m\}.$

\begin{assumption}
\label{chap2strongLC}
 For each $k=1,\hdots,N,$ denote by $u_{S_k}$ the vector with components $u_i$ with $i\in S_k.$ Assume that the strengthened generalized Legendre-Clebsch condition  holds on $\hat{I}_k,$ i.e. 
\benl
-\frac{\partial}{\partial u_{S_k}}\ddot{H}_{u_{S_k}} \succ 0,\quad {\rm  on}\ \Ih_k.
\eenl
\end{assumption}

Hence,  $u_{S_k}$ can be retrieved from equation
\be
\label{chap2ddotHuSk}
\ddot{H}_{u_{S_k}}=0,
\ee
since the latter is affine on $u_{S_k}$
as it has been already pointed out in Section \ref{chap2SecOS}. 
Observe that the expression obtained from \eqref{chap2ddotHuSk} involves only the state variable $\xh$ and the corresponding adjoint state $\ph.$ Hence, it results that $\uh_{S_k}$ is continuous on $\Ih_k$ with finite limits at the endpoints of this interval.
 As the components $\uh_i$ with $i\notin S_k$ are either identically 1 or 0, we conclude that
\be\label{chap2contcont}
\uh\ \text{is continuous on}\ \Ih_k.
\ee

By Assumption \ref{chap2disc} and condition \eqref{chap2contcont} (derived from  Assumption \ref{chap2strongLC}) we get that there exists $\rho \gr 0$ such that 
\be
\label{chap2contpos}
\rho\mi \uh_{i,t} \mi 1-\rho,\quad \text{a.e. on}\ \Ih_k,\ 
\text{for}\ k=1,\hdots,N, \ i\in S_k.
\ee
Next, we present a new control problem obtained in the following way. For each $k=1,\hdots,N,$ we perform the change of time variable that converts the interval $\Ih_k$ into $[0,1]$, afterwards we fix the bang control variables to their bounds and finally, we associate a free control variable to each index in $S_k.$
 More precisely, consider for $k=1,\hdots,N,$  the control variables $u_i^k\in L^{\infty}(0,1;\cR),$ with $i\in S_k,$ and the state variables $x^k\in W^{1,\infty}(0,1;\cR^n).$ 
Let the constants $T_k\in \cR,$ for $k=1,\hdots,N-1,$ which will be considered as state variables of zero-dynamics.
Set $T_0:=0,$ $T_N:=T$ and define the problem on the interval $[0,1]$
\begin{align}
 & \label{chap2cost2} \varphi_0(x^1_0,x^N_1)\rightarrow\min,\\
 & \label{chap2stateeq2} \dot x^k=(T_k-T_{k-1})
\left(\sum_{i\in N_k\cup \{0\}}f_i(x^k)+\sum_{i\in S_k} u_i^{k}f_i(x^k)\right),\quad \ k=1,\hdots,N,\\
 & \label{chap2Teq} \dot T_k=0,\quad \ k=1,\hdots,N-1, \\
 & \label{chap2finalcons2} \eta(x^1_0,x^N_1)=0,\\
 & \label{chap2continuity} x^k_1=x^{k+1}_0,\quad k=1,\hdots,N-1.
\end{align}
Denote by (TP) the problem consisting of equations \eqref{chap2cost2}-\eqref{chap2continuity}.
The link between the original problem (CP) and the transformed one (TP) is given in Lemma \ref{chap2link} below.
Set for each $k=1,\hdots,N:$
\begin{align}
 \label{chap2xk}\xh^k_s &:= \xh(\Th_{k-1}+(\Th_k-\Th_{k-1})s),\quad {\rm for}\ s\in [0,1],\\
 \label{chap2uk}\uh_{i,s}^k &:= \uh_i(\Th_{k-1}+(\Th_k-\Th_{k-1})s),\quad {\rm  for}\ i\in S_k,\ {\rm a.a.}\ s\in [0,1].
\end{align}
Set 
\be\label{chap2Whdef}
\Wh:=((\xh^k)_{k=1}^N,(\uh_ i^k)_{k=1,i\in S_k}^N,(\Th_k)_{k=1}^{N-1}).
\ee
\begin{lemma}\label{chap2link}
If $\wh$ is a Pontryagin minimum of (CP), then $\Wh$  is a weak solution of (TP).
\end{lemma}

\proof
The idea of the proof is to derive the weak optimality of $\Wh$ from the Pontryagin optimality of $\wh$ and condition \eqref{chap2contpos}. Since $\wh$ is a Pontryagin minimum for (CP), there exists $\varepsilon \gr 0$ such that $\wh$ is a minimum in the set of feasible trajectories $w=(x,u)$ satisfying
\be\label{chap2Pontmin}
\|x-\xh\|_{\infty} \mi \varepsilon,\quad
\|u-\uh\|_1 \mi \varepsilon,\quad
\|u-\uh\|_{\infty} \mi 1.
\ee
Consider $\bar\delta,\bar\varepsilon\gr 0,$ and
a feasible solution $((x^k),(u_i^k),(T_k))$ for (TP) such that
\be\label{chap2uikbound}
|T_k-\Th_k|\leq \bar{\delta},\quad \|u_i^k-\uh_i^k\|_{\infty} \mi \bar\varepsilon,\quad \text{for all}\ k=1,\hdots,N.
\ee
We shall relate $\varepsilon$ in \eqref{chap2Pontmin} with $\bar\delta$ and $\bar\varepsilon$ in \eqref{chap2uikbound}. Consider an index $k=1,\hdots,N.$
Denote $I_k:=[T_{k-1},T_k],$ and define for each $i=1,\hdots,m:$
\be\label{chap2udef}
u_{i,t}:= 
\left\{
\ba{cl}
0, &{\rm if}\ t\in I_k\ {\rm  and}\ i\in E_k,\\
u_i^k\left( \frac{t-T_{k-1}}{T_k-T_{k-1}}\right), &{\rm if}\ t\in I_k\ {\rm  and}\ i\in S_k,\\
1, &{\rm if}\ t\in I_k\ {\rm  and}\ i\in N_k.
\ea
\right.
\ee
Let $x$ be the solution of \eqref{chap2stateeq} associated to $u$ and having $x_0=x_0^1.$ We shall prove that $(x,u)$ is feasible for the original problem (CP). Observe that condition 
\eqref{chap2continuity} implies that $x_t=x^k\left(\frac{t-T_{k-1}}{T_k-T_{k-1}}\right)$ when $t\in I_k,$ and thus $x_1=x_1^N.$ It follows that \eqref{chap2finalcons} holds.
We shall check condition \eqref{chap2contcons}.  For $i\in E_k\cup N_k,$ it follows from the definition in \eqref{chap2udef}. Consider now $i\in S_k.$ Since \eqref{chap2contpos} holds, by \eqref{chap2uk} we get
\benl
\rho \mi \uh_{i,s}^k \mi 1-\rho,\quad {\rm  a.e.}\ {\rm  on}\ ]0,1[.
\eenl
Thus, by \eqref{chap2uikbound} and if $\bar\varepsilon \mi \rho,$ we get $0\mi u^k_{i,s}\mi 1$ a.e. on $[0,1].$ This yields
\benl
0\mi u_{i,t} \mi 1,\quad {\rm  a.e.}\ {\rm  on}\ I_k,
\eenl
and thus the feasibility of  $(x,u)$ for (CP).

We now estimate $\|u-\uh\|_1.$ For $k=1,\hdots,N$ and $i\in S_k,$ 
\be\label{chap2uibound}
\ba{rl}
\int_{I_k\cap \hat{I}_k}|u_{i,t}-\uh_{i,t}|\dtt
\leq&
\int_{I_k\cap \hat{I}_k}\  \left|u_i^k\left(\frac{t-T_{k-1}}{T_k-T_{k-1}}\right)-\uh_i^k\left(\frac{t-T_{k-1}}{T_k-T_{k-1}}\right)\right|\dtt
\\
&+
\int_{I_k\cap \hat{I}_k}\  \left|\uh_i^k\left(\frac{t-T_{k-1}}{T_k-T_{k-1}}\right)-\uh_i^k\left(\frac{t-\Th_{k-1}}{\Th_k-\Th_{k-1}}\right)\right|\dtt.
\ea
\ee
Note that, by Assumption \ref{chap2disc} and condition \eqref{chap2contcont}, each $\uh_i^k$ is uniformly continuous on $\Ih_k,$ and thus, there exists $\theta_{ki}\gr0$ such that $|\uh_{i,s}^k-\uh^k_{i,s'}|\mi \bar\varepsilon,$ whenever $|s-s'|\mi \theta_{ki}.$ 
 Set  $\bar\theta:=\min\,\theta_{ki}\gr0.$
Let $\bar\delta$ be such that, if $
|T_k-\Th_k|\mi \bar\delta,
$
then
$
\left| \frac{t-T_{k-1}}{T_k-T_{k-1}} -\frac{t-\Th_{k-1}}{\Th_k-\Th_{k-1}} \right| \mi \bar\theta.
$
From \eqref{chap2uikbound} and \eqref{chap2uibound} we get
\be\label{chap2uibound2}
\int_{I_k\cap \hat{I}_k} |u_{i,t}-\uh_{i,t}|\dtt 
\mi 2\bar\varepsilon\, {\rm meas}\,(I_k\cap \hat{I}_k).
\ee
Assume, w.l.o.g., that $T_k\mi \Th_k$ and note that
\be\label{chap2uibound3}
\int_{T_k}^{\Th_k}|u_{i,t}-\uh_{i,t}|\dtt
\leq 
\int_{T_k}^{\Th_k}\left|u_i^k\left(\frac{t-T_{k-1}}{T_k-T_{k-1}}\right)-\uh_i^k\left(\frac{t-\Th_{k-1}}{\Th_k-\Th_{k-1}}\right)\right|\dtt\mi \bar\delta\,\bar\varepsilon,
\ee
where we used \eqref{chap2uikbound} in the last inequality.
From \eqref{chap2uibound2} and \eqref{chap2uibound3} we get
that $\|u_i-\uh_i\|_1 \mi \bar\varepsilon (2T+(N-1)\bar\delta).$  Thus $\|u-\uh\|_1\mi \varepsilon$ if 
\be\label{chap2epsbound}
\bar\varepsilon(2T+(N-1)\bar\delta) \mi \varepsilon/m.
\ee
We conclude from \eqref{chap2Pontmin} that $((x^k),(u_i^k),(T_k))$ is a minimum on the set of feasible points satisfying  \eqref{chap2uikbound} and \eqref{chap2epsbound}. Thus $\Wh$  is a weak solution of (TP), as it was to be proved.
\qed


We shall next propose a shooting function associated to (TP).
The pre-Hamiltonian of the latter is
\benl
\tilde H:=\sum_{k=1}^N (T_k-T_{k-1})H^k,
\eenl
where, denoting by $p^k$ the costate variable associated to $x^k,$ 
\be\label{chap2Hk}
H^k:=p^k\left(\sum_{i\in N_k\cup \{0\}}f_i(x^k)+\sum_{i\in S_k} u_i^{k}f_i(x^k)\right).
\ee
 Observe that Assumption \ref{chap2strongLC} made on $\uh$ yields
\benl
-\frac{\partial}{\partial u}\ddot{\tilde H}_u \succ 0,\quad {\rm on}\ [0,1],
\eenl
i.e. the strengthened generalized Legendre-Clebsch condition holds in problem (TP) at $\wh.$ Hence we can define the shooting function for (TP) as it was done in Section \ref{chap2SecShoot} for (P).

The endpoint Lagrangian is
\benl 
\tilde\ell:=\varphi_0(x^1_0,x^N_1)+\sum_{j=1}^{d_{\eta}} \beta_j\eta_j(x^1_0,x^N_1)
+\sum_{k=1}^{N-1} \theta_k(x_1^{k}-x_0^{k+1}).
\eenl
The costate equation for $p^k$ is given by
\benl
\dot{p}^k =-(T_k-T_{k-1})D_{x^k}H^k,
\eenl
with endpoint conditions
\begin{gather}
\label{chap2p10}
p_0^1=-D_{x^1_0}\tilde\ell=- D_{x^1_0}\varphi_0
-\sum_{j=1}^{d_{\eta}} \beta_jD_{x^1_0}\eta_j,
\\
\ba{rl}\label{chap2pkboundary}
p^k_1&= \theta^k,\quad {\rm  for}\ k=1,\hdots,N-1,\\
p^k_0&= \theta^{k-1},\quad {\rm  for}\ k=2,\hdots,N,
\ea
\\
\label{chap2pN1}
p_1^N=D_{x^N_1}\tilde\ell= D_{x^N_1}\varphi_0
+\sum_{j=1}^{d_{\eta}} \beta_jD_{x^N_1}\eta_j.
\end{gather}
For the costate variables $p^{T_k}$ associated with $T_k$ we get the equations
\be 
\label{chap2pTk}
\dot{p}^{T_k}=-H^k+H^{k+1},\quad p_0^{T_k}=0,\quad p_1^{T_k}=0,\quad {\rm  for}\ k=1,\hdots,N-1.
\ee
\begin{remarks}
\label{chap2rempTk}
We can sum up the conditions in \eqref{chap2pTk} integrating the first one and obtaining 
$
\int_0^1(H^{k+1}-H^k)\dtt=0,
$
and hence, since $H^k$ is constant on the optimal trajectory, we get the equivalent condition
\be \label{chap2Hcont}
H^k_1=H^{k+1}_0,\quad {\rm  for}\ k=1,\hdots,N-1.
\ee
So we can remove the shooting variable $p^{T_k}$ and  keep the continuity condition on the pre-Hamiltonian.
\end{remarks}

Observe that  \eqref{chap2continuity} and \eqref{chap2pkboundary} imply the continuity of the two functions obtained by concatenating the states and the costates,  i.e. the continuity of $X$ and $P$ defined by
\begin{gather*}
X_0:= x^1_0,\ X_s:=x^k(s-(k-1)),\quad {\rm  for}\ s\in (k-1,k],\ k=1,\hdots,N,
\\
P_0:= p^1_0,\ P_s:=p^k(s-(k-1)),\quad {\rm  for}\ s\in (k-1,k],\ k=1,\hdots,N.
\end{gather*}
Thus, while iterating the shooting method, we can either include the conditions \eqref{chap2continuity} and \eqref{chap2pkboundary} in the definition of the shooting function or integrate the differential equations for $x^k$ and $p^k$ from the values $x^{k-1}_1$ and $p^{k-1}_1$ previously obtained. The latter option reduces the number of variables and hence the size of the problem, but is less stable. We shall present below the shooting function for the more stable case.
To this end define the $n\times n-$matrix 
\benl 
A^k:=\sum_{i\in N_k\cup \{0\}} f_i'(\xh^k) + \sum_{i\in S_k} \uh_i^kf'_i(\xh^k),
\eenl
the $n\times |S_k|-$matrix $B^k$ with columns $f_i(\xh^k)$ with $i\in S_k,$ and 
\benl 
B^k_1:=A^k B^k-\ddt B^k.
\eenl
We shall denote by $g_i(x^k, u^k)$ the $i$th. column of $B^k_1$ for each $i$ in $S_k.$ Here $u^k$ is the $|S_k|-$dimensional vector of components $u_i^k.$
The resulting shooting function for (TP) is given by
\be\label{chap2shoot2}
\ba{rl}
\shoot\colon\cR^{Nn+N-1}\times\cR^{Nn+d_{\eta},*}
&\rightarrow 
\cR^{d_{\eta}+(N-1)n}\times \cR^{(N+1)n+N-1+2\sum |S_k|,*},\\ 
\begin{pmatrix} 
(x^k_0),(T_k),(p^k_0),\beta 
\end{pmatrix}
 =:\nu
&\mapsto
\shoot(\nu):=
\left(
\ba{c}
\vspace{5pt}\eta(x^1_0,x^N_1)\\
\vspace{5pt}(x^k_1-x^{k+1}_0)_{k=1,\hdots,N-1} \\
\vspace{5pt}p^1_0 +D_{x^1_0}\tilde\ell\lam(x^1_0,x^N_1)\\
\vspace{5pt}(p^k_1-p^{k+1}_0)_{k=1,\hdots,N-1}\\
\vspace{5pt}p^N_1 -D_{x^N_1}\tilde\ell\lam(x^1_0,x^N_1)\\
\vspace{5pt} (H_1^k-H_0^{k+1})_{k=1,\hdots,N-1}\\
\vspace{5pt}(p^k_0 f_i(x^k_0))_{k=1,\hdots,N,\ i\in S_k}\\
(p^k_0 g_i(x^k_0,u^k_0))_{k=1,\hdots,N,\ i\in S_k}
\ea
\right).
\ea
\ee
Here, we put both conditions $\tilde H_u=0$ and $\dot{\tilde H}_u=0$ at the beginning of the interval since we have already pointed out in Remark \ref{chap2RemBound} that all the possible choices were equivalent.

Since problem (TP) has the same structure than problem (P) in Section \ref{chap2SecPb}, i.e. they both have free control variable (initial-final constraints), we can apply Theorem \ref{chap2wp} and obtain the analogous result below.

\begin{theorem}\label{chap2cwp}
 Assume that $\wh$ is a Pontryagin minimum of (CP) such that $\Wh$ defined in \eqref{chap2Whdef} satisfies condition \eqref{chap2unifpos} for problem (TP). Then the shooting algorithm for (TP) is locally quadratically convergent.
\end{theorem}

\begin{remarks}
Once system \eqref{chap2shoot2} is obtained, observe that two numerical implementations can be done: one integrating each variable on the interval $[0,1]$ and the other one, going back to the original interval $[0,T],$ 
and using implicitly the continuity conditions \eqref{chap2continuity}, \eqref{chap2pkboundary} and \eqref{chap2Hcont} at each switching time.
The latter implementation is done in the numerical tests of Section \ref{chap2SecNum} below. In this case, the sensitivity with respect to the switching times is obtained from the derivative of the shooting function.
\end{remarks}


\subsection{Reduced Systems}\label{chap2SecRed}

In some cases we can show that some of the conditions imposed to the shooting function in \eqref{chap2shoot2} are redundant. Hence, they can be removed from the formulation yielding a smaller system that we will refer as \textit{reduced system} and which is associated to a \textit{reduced shooting function.} 

Recall that, when defining $\shoot,$ we are implicitly imposing that $\ddot{\tilde H}_u\equiv 0.$ The latter condition together with $\dot{\tilde H}_{u,0}=\tilde H_{u,1}=0,$ both included in the definition of $\shoot,$ imply that 
$\dot{\tilde H}_{u} \equiv \tilde H_u\equiv 0.$ 
Hence,
\be
\label{chap2pk1}
p_1^k  f_i(x_1^k) = p_1^k g_i(x_1^k,u_1^k)=0,\quad \text{for}\ k=1,\hdots,N,
\ i\in S_k,
\ee
and, in view of the continuity conditions \eqref{chap2continuity} and \eqref{chap2pkboundary}, 
\be
\label{chap2pk11}
p_0^{k+1}  f_i(x_0^{k+1}) = p_0^{k+1} g_i(x_0^{k+1},u_0^{k+1})=0,
\quad \text{for}\ k=1,\hdots,N-1,
\ i\in S_k.
\ee
Therefore, if a component of the control is singular on $I_k$ and remains being singular on $I_{k+1},$ then there is no need to impose the boundary conditions on $\tilde H_u$ and $\dot{\tilde H}_{u}$ since they are a consequence of the continuity conditions and the implicit equation $\ddot{\tilde H}_u\equiv 0.$

Observe now that from \eqref{chap2Hk}, \eqref{chap2shoot2} and previous two equations \eqref{chap2pk1} and \eqref{chap2pk11} we obtain,
\benl
H^k_1 = 
p^k_1  \sum_{N_k\cup \{0\}} f_i(x^k_1) =
p^{k+1}_0 \sum_{N_k\cup \{0\} \backslash S_{k+1}} f_i(x^{k+1}_0).
\eenl
On the other hand, 
\benl
H^{k+1}_0 = 
p^{k+1}_0 \sum_{N_{k+1}\cup \{0\} \backslash S_{k}} f_i(x^{k+1}_0).
\eenl
Thus, $H^k_1 = H^{k+1}_0$ if 
$N_k\cup \{0\} \backslash S_{k+1} = N_{k+1}\cup \{0\} \backslash S_{k}.$ The latter equality holds if and only if at instant $T_k$ all the switchings are either bang-to-singular  or singular-to-bang.

\begin{definition}[Reduced shooting function]
We call \textit{reduced shooting function} and we denote it by $\shoot^r$ the function obtained from $\shoot$ defined in \eqref{chap2shoot2} by removing the condition $H^k_1 = H^{k+1}_0$ whenever all the switchings occurring at $T_k$ are either bang-to-singular  or singular-to-bang, and removing  
\benl
p_0^{k}  f_i(x_0^{k}) = 0,\quad p_0^{k} g_i(x_0^{k},u_0^{k})=0,
\eenl
for $k=2,\hdots,N$ and $i\in S_{k-1}\cap S_{k}.$
\end{definition}

\subsection{Square Systems}\label{chap2SSecSquare}

The reduced system above-presented can occasionally result \textit{square}, in the sense that the reduced function $\shoot^r$ has as many variables as outputs. This situation occurs, e.g., in problems 1 and 3 of Section \ref{chap2SecNum}. 
The fact that the reduced system turns out to be square is a consequence of the structure of the optimal solution.
In general, the optimal solution $\uh$ yields a square reduced system if and only if each singular arc is in the interior of $[0,T]$ and at each switching time only one control component switches. This can be interpreted as follows: each singular arc contributes to the formulation with two inputs that are its entry and exit times, and with two outputs that correspond to $p_0^{k}  f_i(x_0^{k}) =g_i(x_0^{k},u_0^{k})=0,$ being $I_k$ the first interval where the component is singular and $i$ the index of the analyzed component. On the other hand, whenever a bang-to-bang transition occurs, it contributes to the formulation with one input for the switching time and one output associated to the continuity of the pre-Hamiltonian (which is sometimes expressed as a zero of the switching function).



\section{Stability under Data Perturbation}\label{chap2SecSta}

In this section, we investigate the stability of the optimal solution under data perturbation. We shall prove that, under condition \eqref{chap2unifpos}, the solution is stable under small perturbations of the data functions $\varphi_0,$ $f_i$ and $\eta.$ 
Assume for this stability analysis that the shooting system of the studied problem can be reduced to a square one. 
We gave a description of this situation in Subsection \ref{chap2SSecSquare}.
Even if the above-mentioned square systems appear in control constrained problems, we start this section by establishing a stability result of the optimal solution for an unconstrained problem. Afterwards, in Subsection \ref{chap2StabCons}, we apply the latter result to problem (TP) and this way we obtain a stability result for the control constrained problem (CP).

\subsection{Unconstrained Control Case}
Consider then problem (P) presented in Section \ref{chap2SecPb}, and the family of problems depending on the real parameter $\mu$ given by:

\be\label{chap2Pmu}\tag{P$_{\mu}$}
\ba{l}
\varphi_0^{\mu}(x_0,x_T)\rightarrow \min,\\
\dot x_t=\ds\sum_{i=0}^m u_{i,t}f_i^{\mu}(x_t),\quad {\rm  a.e.}\ {\rm on}\ [0,T],\\
\eta^{\mu}(x_0,x_T)=0.
\ea
\ee
Assume that $\varphi_0^{\mu}:\cR^{2n+1} \rightarrow \cR$ and $\eta^{\mu}:\cR^{2n+1}\rightarrow \cR^{d_{\eta}}$ 
have Lipschitz-continuous second derivatives in the variable $(x_0,x_T)$ and continuously differentiable with respect to $\mu, $ and $f_i^{\mu}:\cR^{n+1}\rightarrow \cR^n$ is twice continuously differentiable with respect to $x$ and continuously differentiable with respect to the parameter $\mu.$ In this formulation, the problem (P$_0$) associated to $\mu=0$ coincides with (P), i.e. 
$\varphi_0^0=\varphi_0,$  $f_i^0=f_i$ for $i=0,\hdots,m$ and $\eta^0=\eta.$ 
Recall \eqref{chap2unifpos} in Theorem \ref{chap2sc2}, and write the analogous condition for \eqref{chap2Pmu} as follows:
\be\label{chap2unifposmu}
\bar\Omega^{\mu}(\xi,y,h)\geq \rho\gamma(\xi_0,y,h),\quad {\rm  on}\ \P_2^{\mu},
\ee
where $\bar\Omega^{\mu}$ and $\P_2^{\mu}$ are the second variation and critical cone associated to \eqref{chap2Pmu}, respectively.
Let $\shoot^{\mu}$ be the shooting function for (P$_{\mu}$).
Thus, we can write
\benl
\shoot^{\mu}:\cR^M\times \cR \rightarrow \cR^M,\quad
(\,\nu\,,\,\mu\,)\mapsto \shoot^{\mu}(\nu),
\eenl
where we indicate with $M$ the dimension of the domain of $\shoot.$ The following stability result will be established.
\begin{theorem}[Stability of the optimal solution]\label{chap2TheStab}
Assume that the shooting system generated by problem (P) is square and let $\wh$ be a solution satisfying the uniform positivity condition \eqref{chap2unifpos}. Then there exists a neighborhood $\J\subset \cR$ of 0, 
and a continuous differentiable mapping  $\mu\mapsto w^{\mu}=(x^{\mu},u^{\mu}),$ from $\J$ to  $\W,$ where $w^{\mu}$ is a weak solution for \eqref{chap2Pmu}.
Furthermore, $w^{\mu}$ verifies the uniform positivity \eqref{chap2unifposmu}. Therefore, in view of Theorems \ref{chap2sc2} and \ref{chap2wp}, the $\gamma-$ growth holds, and the shooting algorithm for $(P^{\mu})$ is locally quadratically convergent.
\end{theorem}

Let us start showing the following stability result for the family of shooting functions $\{\shoot^{\mu}\}.$
\begin{lemma}\label{chap2LemmaSta1}
Under the hypotheses of Theorem \ref{chap2TheStab}, 
 there exists a neighborhood $\I\subset \cR$ of 0
and a continuous differentiable mapping $\mu\mapsto \nu^{\mu}=(x_0^{\mu},p_0^{\mu},\beta^{\mu}),$ from $\I$ to  $\cR^M,$ such that $\shoot^{\mu}(\nu^{\mu})=0.$
Furthermore, the solutions $(x^{\mu},u^{\mu},p^{\mu})$ of the system of equations \eqref{chap2stateeq}, \eqref{chap2costateeq}, \eqref{chap2ddotPhi} with initial condition $(x_0^{\mu},p_0^{\mu})$ and associated multiplier $\beta^{\mu}$ provide a family of feasible trajectories $w^{\mu}:=(x^{\mu},u^{\mu})$ verifying
\be\label{chap2estSta}
\|x^{\mu}-\xh\|_{\infty} + \|u^{\mu}-\uh\|_{\infty}+ \|p^{\mu}-\ph\|_{\infty} + |\beta^{\mu}-\hat\beta|
= \mathcal{O}(\mu).
\ee
\end{lemma}

\proof
Since \eqref{chap2unifpos} holds, the result in Theorem \ref{chap2wp} yields the non-singularity of the square matrix $D_{\nu}S^0(\hat\nu).$ Hence, the Implicit Function Theorem is applicable and we can then guarantee the existence of a neighborhood $\B\subset \cR^M$ of $\hat\nu,$ a neighborhood $\I\subset \cR$ of 0, and a continuously differentiable function $\Gamma:\I\rightarrow \B$ such that
\be \label{chap2shootGamma}
\shoot^{\mu}(\Gamma(\mu))=0,\quad {\rm  for}\ {\rm  all}\ \mu\in\I.
\ee
Finally, write $\nu^{\mu}:=\Gamma(\mu)$ and use the continuity of $D\Gamma$ on $\I$ to get the first part of the statement.

The feasibility of $w^{\mu}$ holds since equation \eqref{chap2shootGamma} is verified. Finally, the estimation \eqref{chap2estSta} follows from the stability of the system of
differential equation provided by the shooting method.
\qed

Once we obtained the existence of this  $w^{\mu}$  feasible  for $(P^{\mu}),$ we may wonder whether it is locally optimal. For this aim, we shall investigate the stability of the sufficient condition \eqref{chap2unifpos}. Denote by $\bar\Omega^{\mu}$ and $\P_2^{\mu}$ the quadratic mapping and critical cone related to $(P^{\mu}),$ respectively. Given that all the functions involved in $\bar\Omega^{\mu}$ are continuously differentiable with respect to $\mu,$ the mapping $\bar\Omega^{\mu}$ itself is continuously differentiable with respect to $\mu.$ For the perturbed cone we get the following approximation result.
\begin{lemma}
\label{chap2LemmaSta2}
 Assume the same hypotheses as in Theorem \ref{chap2TheStab}.
Take $\mu\in\I$ and \\
$(\xi^{\mu},y^{\mu},h^{\mu})\in\P_2^{\mu}.$ Then there exists
$(\xi,y,h)\in\P_2$ such that
\benl
|\xi_0^{\mu}-\xi_0| + 
\|y^{\mu}-y\|_{2} +
|h^{\mu}-h|
=\mathcal{O}(\mu).
\eenl
\end{lemma}
The definition below will be useful in the proof of previous Lemma.
\begin{definition}
 Define the function
$
\bar\eta:\U\times\cR^n\rightarrow \cR^{d_{\eta}},
$ 
given by
\benl
\bar\eta(u,x_0):=\eta(x_0,x_T),
\eenl
where $x$ is the solution of \eqref{chap2stateeq} associated to $(u,x_0).$
\end{definition}

\proof
\ [of Lemma \ref{chap2LemmaSta2}]
Recall that $D\bar\eta(\uh,\xh_0)$ is onto by Assumption \ref{chap2lambdaunique}.
Call back the definition of the critical cone $\C$ given in \eqref{chap2coneC}, and note that we can rewrite it as
$
\C=\{(z,v)\in\W:\, \G(z,v)=0\}=\Ker \G,
$
where $\G$ is an onto linear application from $\W$ to $\cR^{d_{\eta}},$ defined by $\G(z,v):=D\eta(\xh_0,\xh_T)(z_0,z_T).$
In view of Goh's Transformation \eqref{chap2ydef},
\benl
D\eta(\xh_0,\xh_T)(z_0,z_T)=D\eta(\xh_0,\xh_T)(\xi_0,\xi_T+B_Ty_T),
\eenl
for $(z,v)\in \W$ and $(\xi,y)$ being its corresponding transformed direction.
Thus, the cone $\P_2$ can be written as
$
\P_2=\{\zeta\in \H:\, \K(\zeta)=0\}=\Ker \K,
$
with $\zeta:=(\xi,y,h),$ $\H:=\X_2\times \U_2\times \cR^n,$ 
and $\K(\zeta):=D\eta(\xh_0,\xh_T)(\xi_0,\xi_T+B_Th).$ Then $\K\in \L(\H,\cR^{d_{\eta}})$ and it is surjective. Analogously,
one has the identity $
\P_2^{\mu}=\{\zeta\in \H:\, \K^{\mu}(\zeta)=0\}=\Ker \K^{\mu},
$
with
\be
\label{chap2Kmu}
\|\K^{\mu}-\K\|_{\L(\H,\cR^{d_{\eta}})}=\mathcal{O}(\mu).
\ee
Let us now prove the desired stability property. Take $\zeta^{\mu}\in \P_2^{\mu}=\Ker \K^{\mu}$ having $\|\zeta\|_{\H}^{\mu}=1.$ Hence
$
\label{chap2Kzeta}
\K(\zeta^{\mu})=\K^{\mu}(\zeta^{\mu})+(\K-\K^{\mu})(\zeta^{\mu}),
$
and by estimation \eqref{chap2Kmu},
\be
\label{chap2Kzetanorm}
|\K(\zeta^{\mu})|=\mathcal{O}(\mu).
\ee
Observe that, since $\H=\Ker \K \oplus {\rm Im}\, \K\tras,$ there exists $\zeta^{\mu,*}\in \H^*$ such that
\be
\label{chap2zeta}
\zeta:=\zeta^{\mu} + \K\tras(\zeta^{\mu,*})\in \Ker\K.
\ee
This yields
$
0=\K(\zeta)=\K(\zeta^{\mu})+\K\K\tras (\zeta^{\mu,*})
=
(\K-\K^{\mu})(\zeta^{\mu})+\K\K\tras (\zeta^{\mu,*}).
$
Given that $\K$ is onto, the operator $\K\K\tras$ is invertible and thus
\benl
\zeta^{\mu,*}=-(\K\K\tras)^{-1}(\K-\K^{\mu})(\zeta^{\mu}).
\eenl
The estimation \eqref{chap2Kzetanorm} above implies
$
\|\zeta^{\mu,*}\|_{\H^*}=\mathcal{O}(\mu).
$
It follows then from \eqref{chap2zeta} that 
$
\|\zeta^{\mu}-\zeta\|_{\H}=\mathcal{O}(\mu),
$ and therefore, the desired result holds.
\qed

\proof
\ [of Theorem \ref{chap2TheStab}]
We shall begin by observing that Lemma \ref{chap2LemmaSta1} provides a neighborhood $\I$ and a class of solutions $\{(x^{\mu},u^{\mu},p^{\mu},\beta^{\mu})\}_{\mu\in\I}$ satisfying \eqref{chap2estSta}.
We shall prove  that $w^{\mu}=(x^{\mu},u^{\mu})$ satisfies the sufficient condition \eqref{chap2unifposmu} close to 0.

Suppose on the contrary that there exists a sequence of parameters $\mu_k\rightarrow 0$ and critical directions  $(\xi^{\mu_k},y^{\mu_k},h^{\mu_k})\in\P_2^{\mu_k}$ with $\gamma(\xi_0^{\mu_k},y^{\mu_k},h^{\mu_k})=1,$ such that
\benl
\bar\Omega^{\mu_k}(\xi^{\mu_k},y^{\mu_k},h^{\mu_k})\leq o(1).
\eenl
Since $\bar\Omega^{\mu}$ is Lipschitz-continuous in $\mu,$ from previous inequality we get
\be
\label{chap2Omegabarmuk}
\bar\Omega (\xi^{\mu_k},y^{\mu_k},h^{\mu_k})\leq o(1).
\ee
In view of Lemma \ref{chap2LemmaSta2}, there exists for each $k,$ a direction $(\xi^k,y^k,h^k)\in\P_2$ satisfying
\be\label{chap2ok}
|\xi_0^k-\xi_0^{\mu_k}|+\|y^k-y^{\mu_k}\|_2+|h^k-h^{\mu_k}|=\mathcal{O}(\mu_k).
\ee
Hence, by inequality \eqref{chap2Omegabarmuk} and given that $\wh$ satisfies \eqref{chap2unifpos},
\benl
\rho\gamma(\xi_0^k,y^k,h^k)\leq \bar\Omega(\xi^k,y^k,h^k)\leq o(1).
\eenl
However, the left hand-side of last inequality cannot go to 0 since $(\xi_0^k,y^k,h^k)$ is close to $(\xi_0^{\mu_k},y^{\mu_k},h^{\mu_k})$ by estimation \eqref{chap2ok}, and the elements of the latter sequence have unit norm. This leads to a contradiction. Hence, the result follows.
\qed


\subsection{Control Constrained Case}\label{chap2StabCons}

In this paragraph, we aim to investigate the stability of the shooting algorithm applied to the problem with control bounds (CP) studied in Section \ref{chap2SecCons}.
Observe that previous Theorem \ref{chap2TheStab} guarantees the weak optimality for the perturbed problem when the control constraints are absent. In case we have control constraints, this stability result is applied to the transformed problem (TP) (given by equations \eqref{chap2cost2}-\eqref{chap2continuity} of Section \ref{chap2SecCons}) yielding a similar stability property,  but for which the nominal point and the perturbed ones are  weak optimal for (TP).  This means that they are optimal in the class of trajectories having the same control structure, and switching times and singular arcs sufficiently close in $L^{\infty}.$ 
A trajectory satisfying optimality in this sense will be called \textit{weak-structural optimal,} and a formal definition would be as follows.
\begin{definition}[Weak-structural optimality]
A feasible trajectory $\wh$ for problem (CP) is called a {\em weak-structural solution} iff its transformed trajectory $\Wh$ given by \eqref{chap2xk}-\eqref{chap2Whdef} is a weak solution of (TP). 
\end{definition}

\begin{theorem}[Sufficient condition for the extended weak minimum in the control constrained case]
 Let $\wh$ be a feasible solution for (CP) satisfying Assumptions \ref{chap2geohyp} and \ref{chap2disc}. Consider the transformed problem (TP) and the corresponding transformed solution $\Wh$ given by \eqref{chap2xk}-\eqref{chap2Whdef}. If $\wh$ satisfies \eqref{chap2unifpos} for (TP), then $\wh$ is an extended weak solution for (CP).
\end{theorem}

\proof
It follows from the sufficient condition in Theorem \ref{chap2sc2} applied to (TP).
\qed

Consider the family of perturbed problems 
\be\label{chap2CPmu}\tag{CP$_{\mu}$}
\ba{l}
\varphi_0^{\mu}(x_0,x_T)\rightarrow \min,\\
\dot x_t=\ds\sum_{i=0}^m u_{i,t}f_i^{\mu}(x_t),\quad {\rm  a.e.}\ {\rm on}  \ [0,T],\\
\eta^{\mu}(x_0,x_T)=0,\\
0\leq u_t\leq 1,\quad {\rm  a.e}\ {\rm  on}\ [0,T].
\ea
\ee
The following stability result follows from Theorem \ref{chap2TheStab}.

\begin{theorem}[Stability in the control constrained case]
\label{chap2TheStabCons}
Assume that the reduced shooting system generated by the problem (CP) is square. Let $\wh$ be a solution of (CP) and  $\{\Th_k\}_{k=1}^N$ its switching times.
Denote by $\Wh$ its transformation via equation \eqref{chap2Whdef}. 
Suppose that $\Wh$ satisfies the uniform positivity condition \eqref{chap2unifpos} for problem (TP). Then there exists a neighborhood $\J\subset \cR$ of 0, such that for every parameter
$\mu\in\J,$ there exists a weak-structural optimal trajectory $w^{\mu}$ of $(CP^{\mu})$ with switching times $\{T^{\mu}_k\}_{k=1}^N,$ satisfying the estimation
\benl
\sum_{k=1}^N|T^{\mu}_k-\Th_k|
+
\sum_{k=1}^N \sum_{i\in S_k} \|u^{\mu}_i-\uh_i\|_{\infty,I^{\mu}_k\cap \Ih_k}
+
\|x^{\mu}-\xh\|_{\infty} = \mathcal{O}(\mu),
\eenl
where $I^{\mu}_k:=[T_{k-1}^{\mu},T_k].$
Furthermore, the transformed perturbed solution $W^{\mu}$ verifies the uniform positivity \eqref{chap2unifposmu} and hence, the quadratic growth in the weak sense for problem (TP) holds, and the shooting algorithm for \eqref{chap2CPmu} is locally quadratically convergent.
\end{theorem}

\subsection{Additional Analysis for the Scalar Control Case}

Consider a particular case where the control $\uh$ is scalar.
The lemma below shows that the perturbed solutions are Pontryagin minima for \eqref{chap2CPmu}, provided that the following assumption holds.
\begin{assumption}
\label{chap2sc}
(a) The switching function $H_u$ is never zero in the interior of a bang arc. Hence, if $\uh=1$ on $]t_1,t_2[\subset [0,T],$ then $H_u\mi 0$ on $]t_1,t_2[,$ and if $\uh=0$ on $]t_1,t_2[,$ then $H_u\gr 0$ on $]t_1,t_2[.$  

(b) If $\Th_k$ is a bang-to-bang switching time then $\dot{H}_u(\Th_k)\neq 0.$
\end{assumption}

The property (a) is called \textit{strict complementarity for the control constraint.}

\begin{lemma}
 Suppose that $\uh$ satisfies Assumption \ref{chap2sc}.
Let $w^{\mu}$ be as in Theorem \ref{chap2TheStabCons} above.
Then $w^{\mu}$ is a Pontryagin minimum for \eqref{chap2CPmu}.
\end{lemma}

\proof
 We intend to prove that $w^{\mu}$ satisfies the minimum condition \eqref{chap2maxcond} given by the Pontryagin Maximum Principle. 
Observe that on the singular arcs, $H_u^{\mu}=0$ since $w^{\mu}$ is the solution associated to a zero of the shooting function.
It suffices then to study the stability of the sign of  $H_u^{\mu}$ on the bang arcs around a switching time.
First suppose that $\uh$ has a bang-to-singular switching at $\Th_k.$ Assume, without any loss of generality, that $\uh\equiv 1$ on $\Ih_k$ and $\uh$ is singular on $[\Th_{k},\Th_{k+1}].$ 
Let us write 
\be
\label{chap2ddotHumu}
\ddot{H}_u^{\mu}=a^{\mu}+u^{\mu}b^{\mu},
\ee
where $a^{\mu}$ and $b^{\mu}:=\frac{\partial}{\partial u} \ddot{H}_u^{\mu}$ are continuous functions on $[0,T],$ and continuously differentiable with respect to $\mu,$ since they
depend on $x^{\mu}$ and $p^{\mu}.$
Assumption \ref{chap2strongLC} yields $b^0\mi 0$ on $[\Th_{k},\Th_{k+1}]$ and, therefore,
\be 
\label{chap2Humu}
b^{\mu}\mi 0,\quad {\rm  on}\ [T^{\mu}_{k},T_{k+1}^{\mu}].
\ee
Due to \eqref{chap2ddotHumu}, the sign of $\ddot{H}_u^{\mu}$ around $T^{\mu}_{k}$ depends on $u^{\mu}(T^{\mu}_{k}+)-u^{\mu}(T^{\mu}_{k}-).$ However, this quantity is negative since $u^{\mu}$ passes from its upper bound to a singular arc. 
From the latter assertion and \eqref{chap2Humu} it follows
\benl
\ddot{H}_u^{\mu}(T^{\mu}_{k}-)\mi 0,
\eenl
 and thus, $H_u^{\mu}$ is concave at the junction time $T^{\mu}_{k}.$ Since $H_u^{\mu}$ is null on $[T^{\mu}_{k},T_{k+1}^{\mu}],$ its concavity implies that it has to be negative before entering this arc. Hence, $w^{\mu}$ respects the minimum condition on the interval $\Ih_k.$

Consider now the case when $\uh$ has a bang-to-bang switching at $\Th_k.$ Let us begin by showing that $H_u^{\mu}(T_k^{\mu})=0.$ Suppose, on the contrary, that we have $H_u^{\mu}(T_k^{\mu})\neq0.$ Then 
$H^{\mu}(T_k^{\mu}+)-H^{\mu}(T_k^{\mu}-)\neq 0,$ contradicting the continuity condition imposed on $H$ in the shooting system.
Hence $H_u^{\mu}(T_k^{\mu})=0.$
On the other hand, since $\dot{H}_u(\Th_k)\neq0$ by Assumption \ref{chap2sc}, the value $\dot{H}_u^{\mu}(T_k^{\mu})$ has the same sign for small values of $\mu.$ This implies that $H_u^{\mu}$ has the same sign that $H_u,$ before and after $T_k^{\mu}$   (or before and after $\Th_k$).
The result follows. 
\qed

\begin{remarks}
 We end this analysis by mentioning that, if the transformed solution $\Wh$ satisfies the uniform positivity \eqref{chap2unifpos} for (TP), then $\wh$ verifies the sufficient condition established in Aronna et al. \cite{ABDL11} and hence it is actually a Pontryagin minimum. This follows from the fact that in condition \eqref{chap2unifpos} we are allowed 
to perturb the switching times, and hence \eqref{chap2unifpos} is more restrictive (or demanding) than the condition in \cite{ABDL11}.
\end{remarks}


\section{Numerical Simulations}\label{chap2SecNum}

Now we aim to check numerically the extended shooting method described above.
More precisely, we want to compare the classical $n \times n$ shooting formulation to an extended formulation with the additional conditions on the pre-Hamiltonian continuity.
We test three problems with singular arcs: a fishing and a regulator problem, and the well-known Goddard problem, which we have already studied in \cite{GeMa06,BonMarTre08}. 
For each problem, we perform a batch of shootings on a large grid around the solution.
We then check the convergence and the solution found, as well as the singular values and condition number of the Jacobian matrix of the shooting function.

\subsection{Test Problems}
\subsubsection{Fishing Problem}
The first example we consider is a fishing problem described in \cite{Clark}.
The state $x_t \in \cR$ represents the fish population (halibut), the control $u_t \in \cR$ is the fishing activity, and the objective is to maximize the net revenue of fishing over a fixed time interval.
The coefficient $(E - {c}/{x})$ takes into account the greater fishing cost for a low fish population. The problem is
\be\tag{P$_1$}
\left \lbrace
\ba{rl}
&\max \ds\int_0^T\,\left(E-{c}/{x_t}\right)\, u_t\, U_{{\rm max}}\dtt,\\
&\dot{x}_t = r\,x_t\,\left(1-{x_t}/{k}\right)\, -\, u_t\, U_{{\rm max}},\\
&0 \le u_t \le 1,\quad {\rm a.e.}\ {\rm on}\ [0,T],\\
&x_0 = 70, \quad x_T\ {\rm  free},\\ 
\ea
\right .
\ee
with $T=10$, $E=1$, $c=17.5$, $r=0.71$, $k=80.5$ and $U_{{\rm max}}=20$.

\begin{remarks}
{The state and control were rescaled by a factor $10^6$ compared to the original data for a better numerical behavior.}
\end{remarks}
\begin{remarks}
Since we have an integral cost, we add a state variable to adapt (P$_1$) to the initial-final cost formulation. It is well-known that its corresponding costate variable is constantly equal to 1. 
\end{remarks}

\begin{figure}
\centering
\begin{center}
\subfigure{\includegraphics [width=6.8cm]{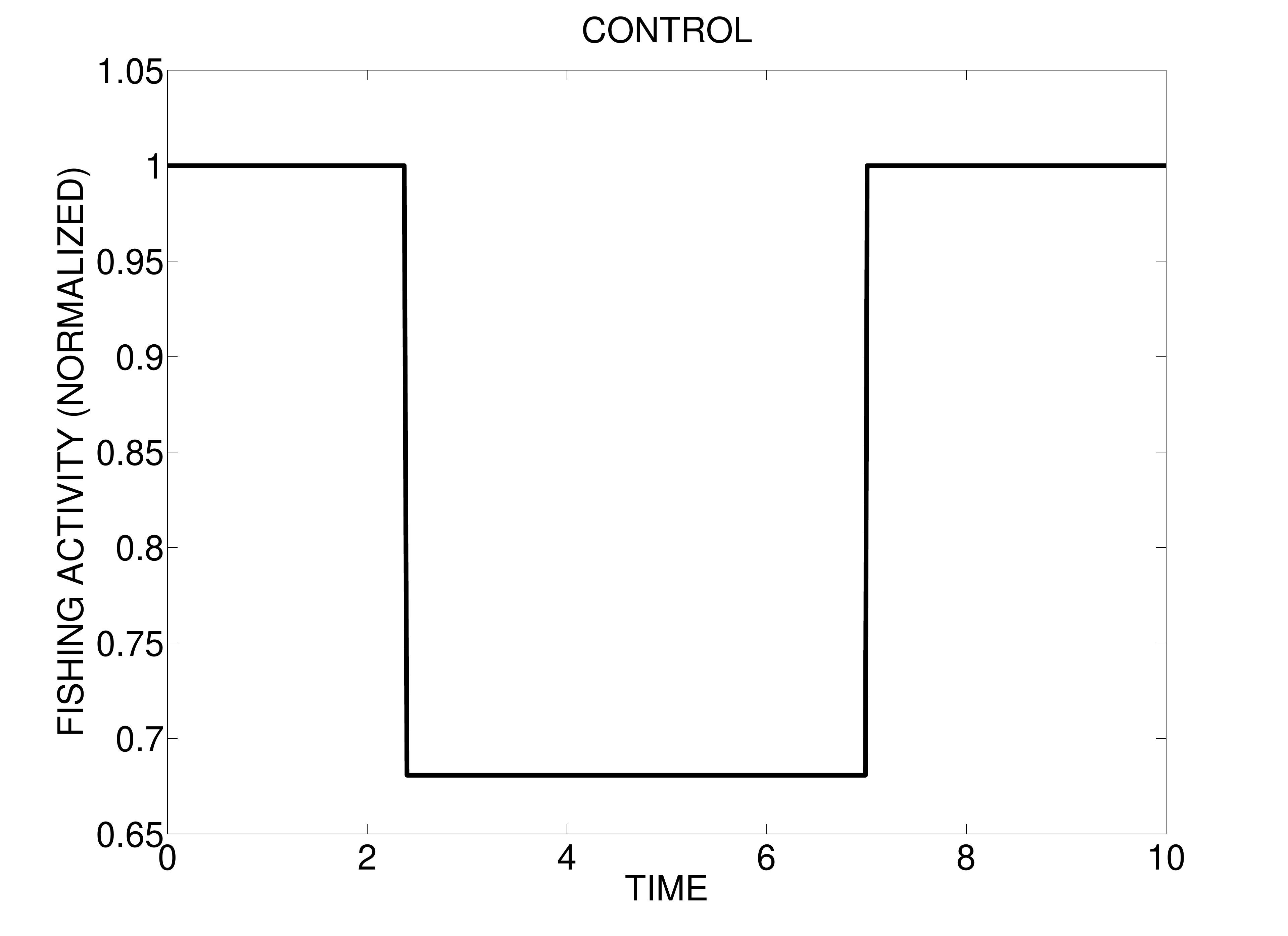}} \\
\subfigure{\includegraphics [width=5.8cm]{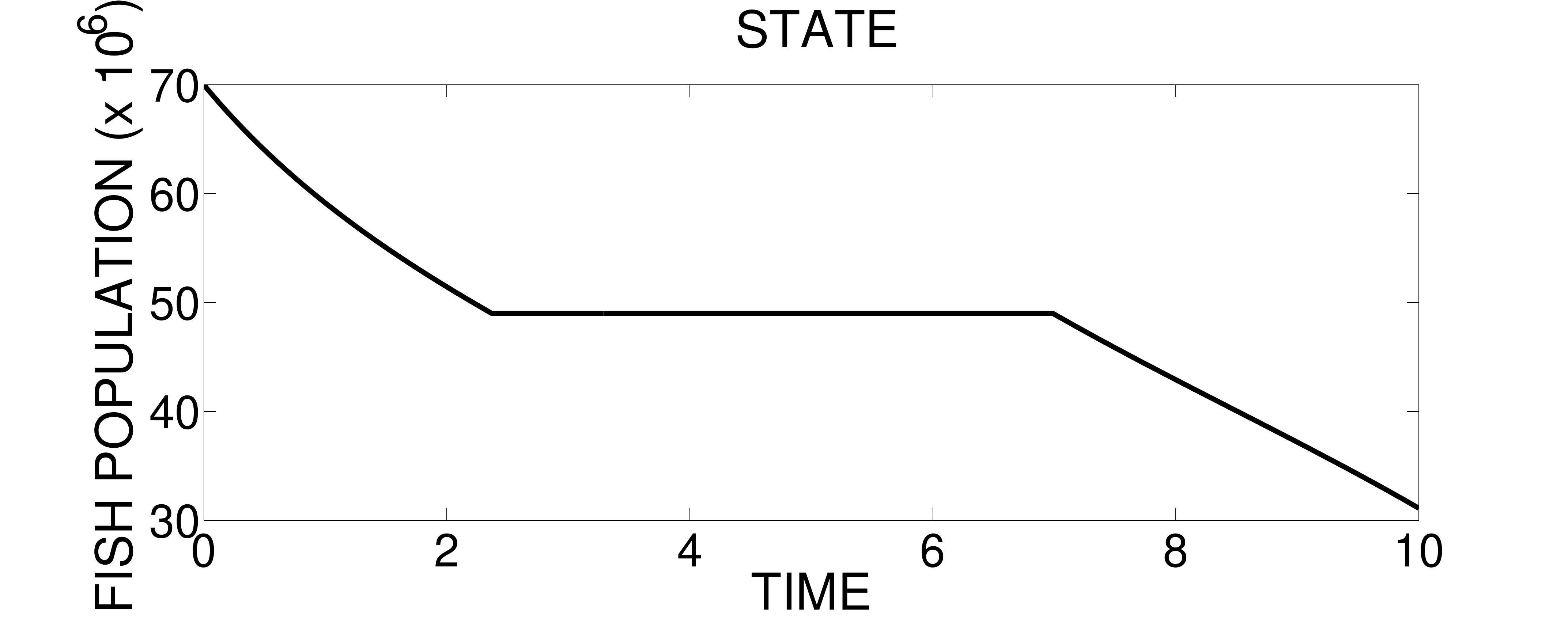}}
\subfigure{\includegraphics [width=5.9cm]{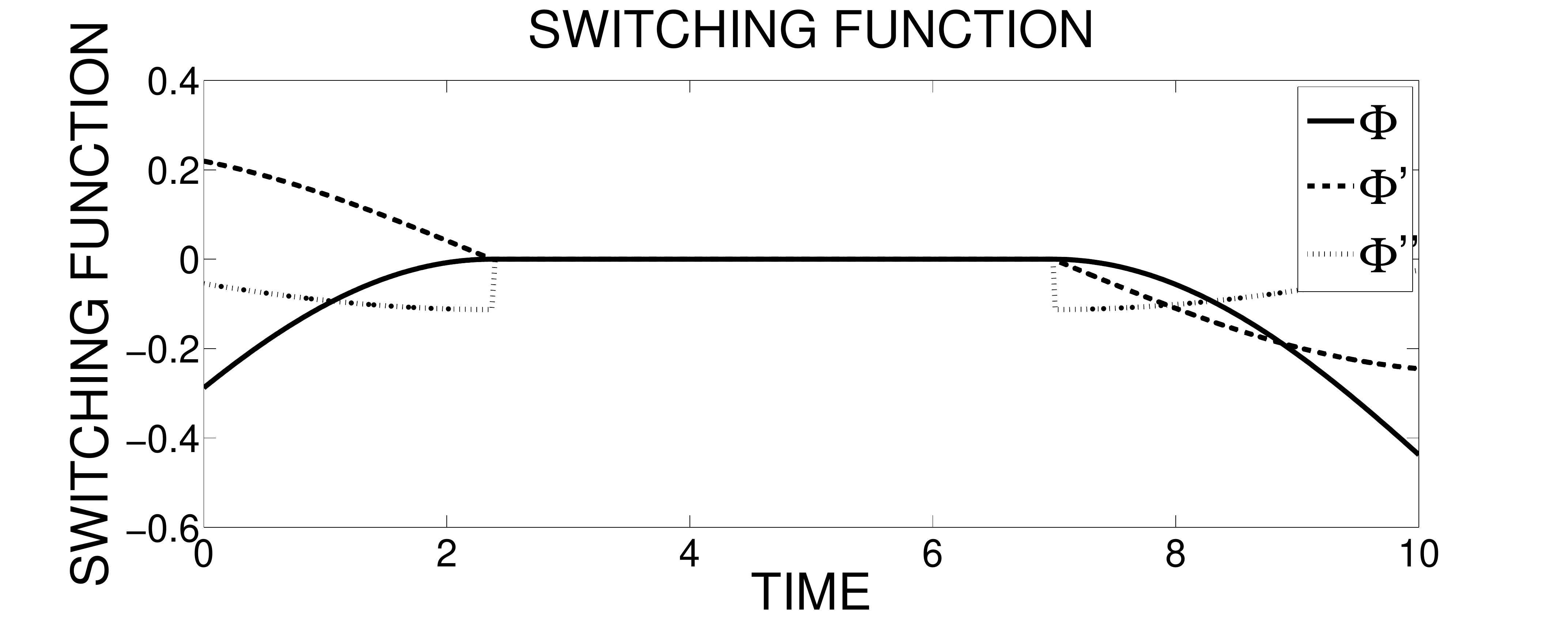}}
\end{center}
\caption{Fishing Problem}
\label{chap2FigFishing}
\end{figure}

The pre-Hamiltonian for this problem is
\benl
H:= (c/x-E)\,u\,U_{{\rm max}}+p[r\,x\,(1-x/k)-u\,U_{{\rm max}}],
\eenl
and hence the switching function
\benl
\Phi_t = D_uH_t = U_{\rm max}({c}/{x_t} - E - p_t),\quad \forall t \in [0,T].
\eenl
The optimal control follows the bang-bang law
\benl
\left \lbrace
\begin{array}{ll}
\uh_t = 0 \quad &{\rm if}\ \Phi_t > 0,\\
\uh_t = 1 \quad &{\rm if}\ \Phi_t < 0.
\end{array}
\right .
\eenl
Over a singular arc, where $\Phi=0$, we assume that the relation $\ddot{\Phi}=0$ gives the expression of the singular control \emph{(t is omitted for clarity)}
\benl
\uh_{{\rm singular}} = \frac{k\ r}{2({c}/{\xh}-\ph)U_{max}} \left(\frac{c}{\xh} - \frac{c}{k} - \ph + \frac{2\ph\xh}{k} - \frac{2\ph\xh^2}{k^2}\right).
\eenl
The solution obtained for (P$_1$) has the structure \textbf{bang-singular-bang}, as shown in Figure \ref{chap2FigFishing}. All the graphics in this article have been done with Matlab.

\noindent \textbf{Shooting formulations.} Assuming the control structure, the shooting unknowns are the initial costate and the endpoints of the singular arc,
\benl
\nu := (p_0, t_1, t_2) \in \cR^3.
\eenl
The classical shooting formulation uses the entry conditions on $t_1$
\benl
\shoot_1(\nu) := (p_T, \Phi_{t_1}, \dot \Phi_{t_1}).
\eenl
The equation $S_1(\nu)=0$ is a square nonlinear system, for which a quasi-Newton method can be used.
Note that, even if there is no explicit condition on $t_2$ in $S$, the value of $p_T$ does depend on $t_2$ via the control switch.

The extended shooting formulation adds two conditions corresponding to the continuity of the pre-Hamiltonian at the junctions between bang and singular arcs.
We denote $[H]_t := H_{t+} - H_{t-}$ the pre-Hamiltonian jump, and define 
\benl
\tilde{\shoot}_1(\nu) :=(p_{10}, \Phi_{t_1}, \dot \Phi_{t_1}, [H]_{t_1}, [H]_{t_2}).
\eenl
To solve $\tilde{S}_1(\nu)=0$ we use a nonlinear least-square algorithm (see paragraph \ref{chap2results} below for more details).


\subsubsection{Regulator Problem} 

\begin{figure}[h]
\centering
\begin{center}
\subfigure{\includegraphics [width=5.8cm]{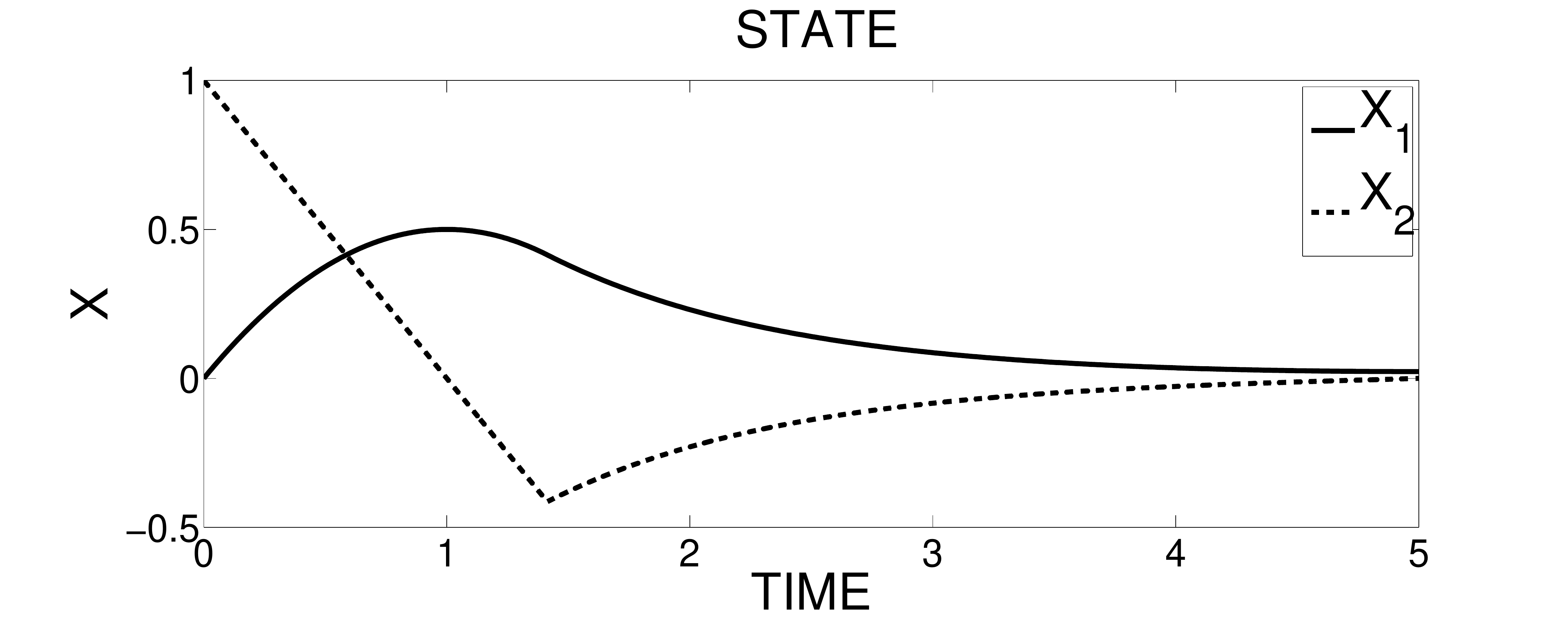}}
\,
\subfigure{\includegraphics [width=5.8cm]{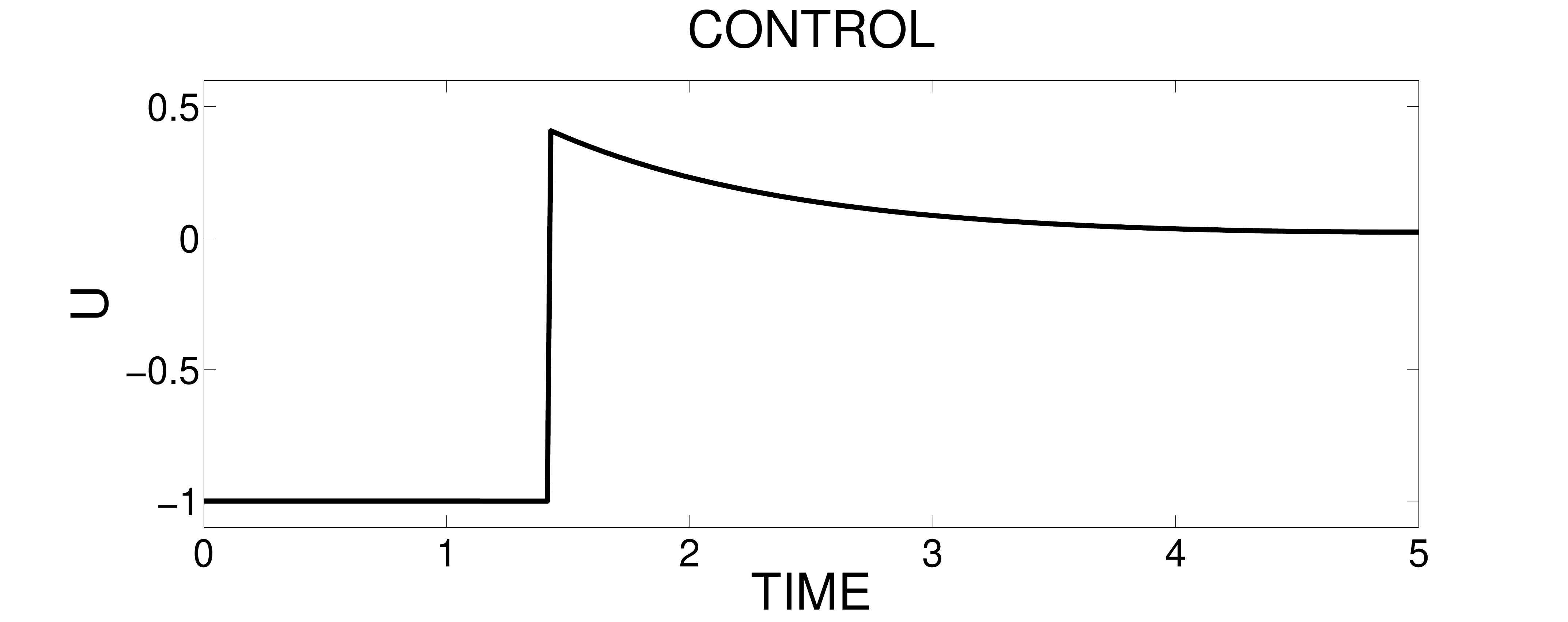}}
\\
\subfigure{\includegraphics [width=7cm]{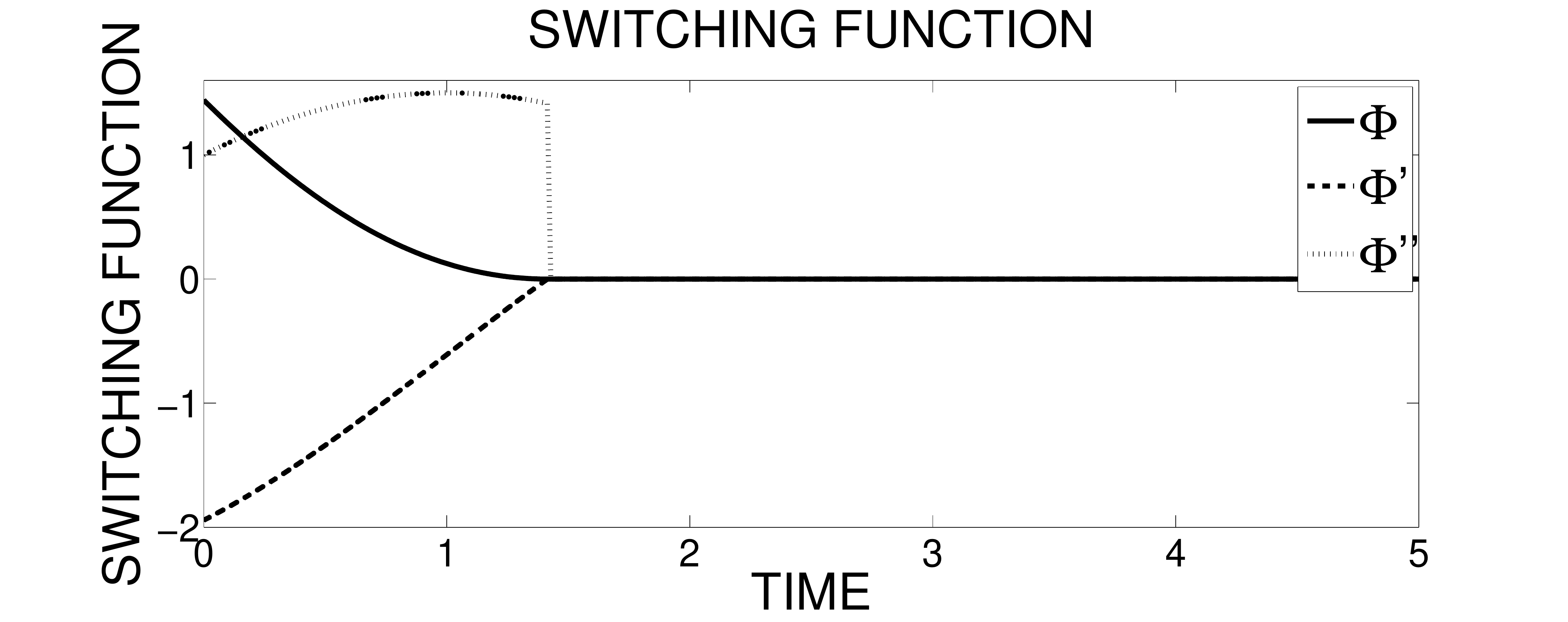}}
\end{center}
\caption{Regulator Problem}
\label{chap2FigRegulator}
\end{figure}

The second example is the linear-quadratic regulator problem described in Aly \cite{Aly78}.
We want to minimize the integral of the sum of the squares of the position and speed of a mobile over a fixed time interval, the control being the acceleration. 
\be \tag{P$_2$}
\left\{
\ba{rl}
&\min \half \ds\int_0^{T}\left( x_{1,t}^2+ x_{2,t}^2  \right) \dtt,\\
&\dot{x}_{1,t} = x_{2,t},\\
&\dot{x}_{2,t} = u_t,\\
&-1 \le u_t \le 1,\quad {\rm  a.e.}\ {\rm  on}\ [0,T],\\
&x_0 = (0,1), \quad x_T\ {\rm  free},\quad T=5.
\ea
\right.
\ee
The corresponding pre-Hamiltonian and the switching function are
\benl
H:= \half( x_1^2+ x_2^2 )  +  p_1x_2+ p_2u,
\eenl
\benl
\Phi_t := D_u H_t = p_{2,t}.
\eenl
The bang-bang optimal control satisfies
\benl
\uh_t = -{\rm sign}\ \ph_{2,t} \quad {\rm if}\ \Phi_t \ne 0.
\eenl
The singular control is again obtained from $\ddot{\Phi}=0$ and verifies
\benl
\uh_{{\rm singular},t} = \xh_{1,t}.
\eenl
The solution for this problem has the structure \textbf{bang-singular}, as shown on Figure \ref{chap2FigRegulator}.

\noindent \textbf{Shooting formulations.} Assuming the control structure, the shooting unknown is
\benl
\nu: = (p_{1,0}, p_{2,0}, t_1) \in \cR^3.
\eenl
For the classical shooting formulation, in order to have a square system, we can, for instance, combine the two entry conditions on $\Phi$ and $\dot \Phi$, since we only have one additional unknown which is the entry time $t_1$. 
Thus we define
\benl
\shoot_2(\nu) := (p_{1,T}, p_{2,T}, \Phi_{t_1}^2 + \dot \Phi_{t_2}^2).
\eenl
The extended formulation does not require such a trick, we simply have
\benl
\tilde{\shoot}_2(\nu) := (p_{1,T}, p_{2,T}, \Phi_{t_1}, \dot \Phi_{t_1}, [H]_{t_1}).
\eenl


\subsubsection{Goddard Problem}

\begin{figure}[!h]
\centering
\begin{center}
\subfigure{\includegraphics [width=5.8cm]{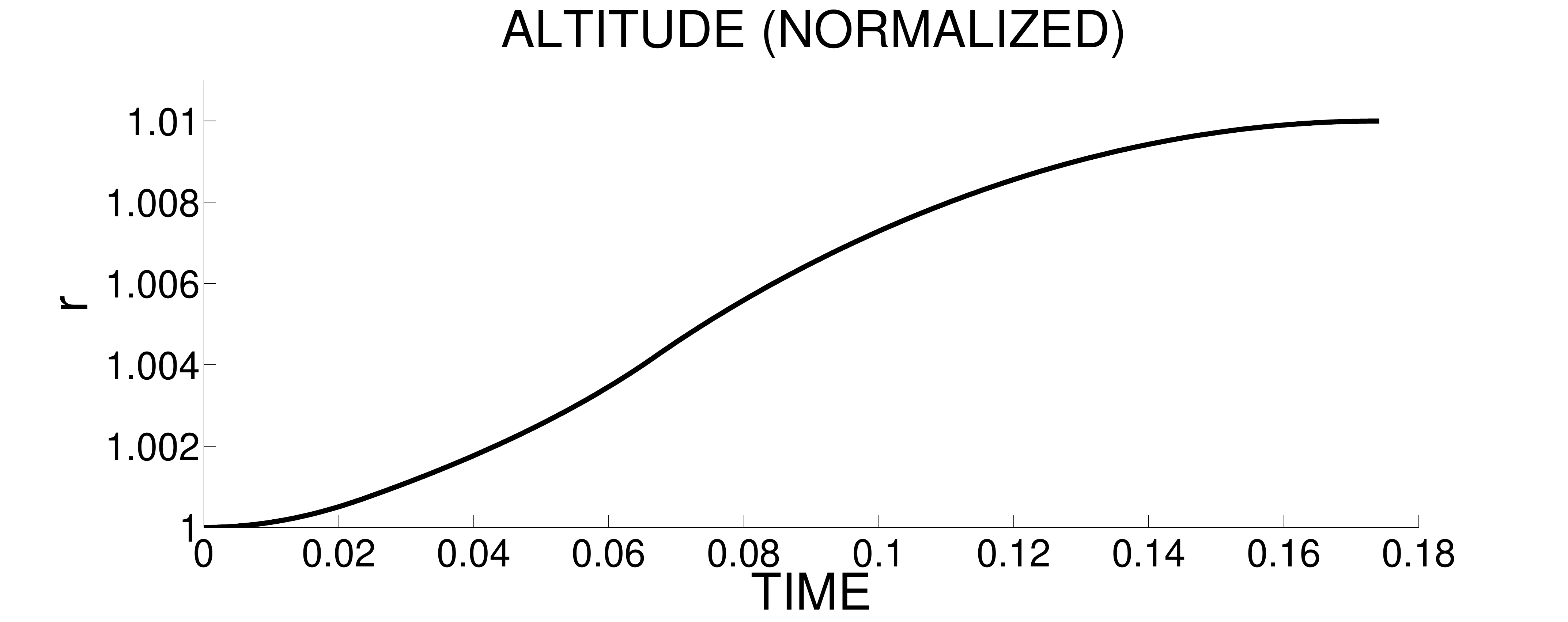}}
\ 
\subfigure{\includegraphics [width=5.8cm]{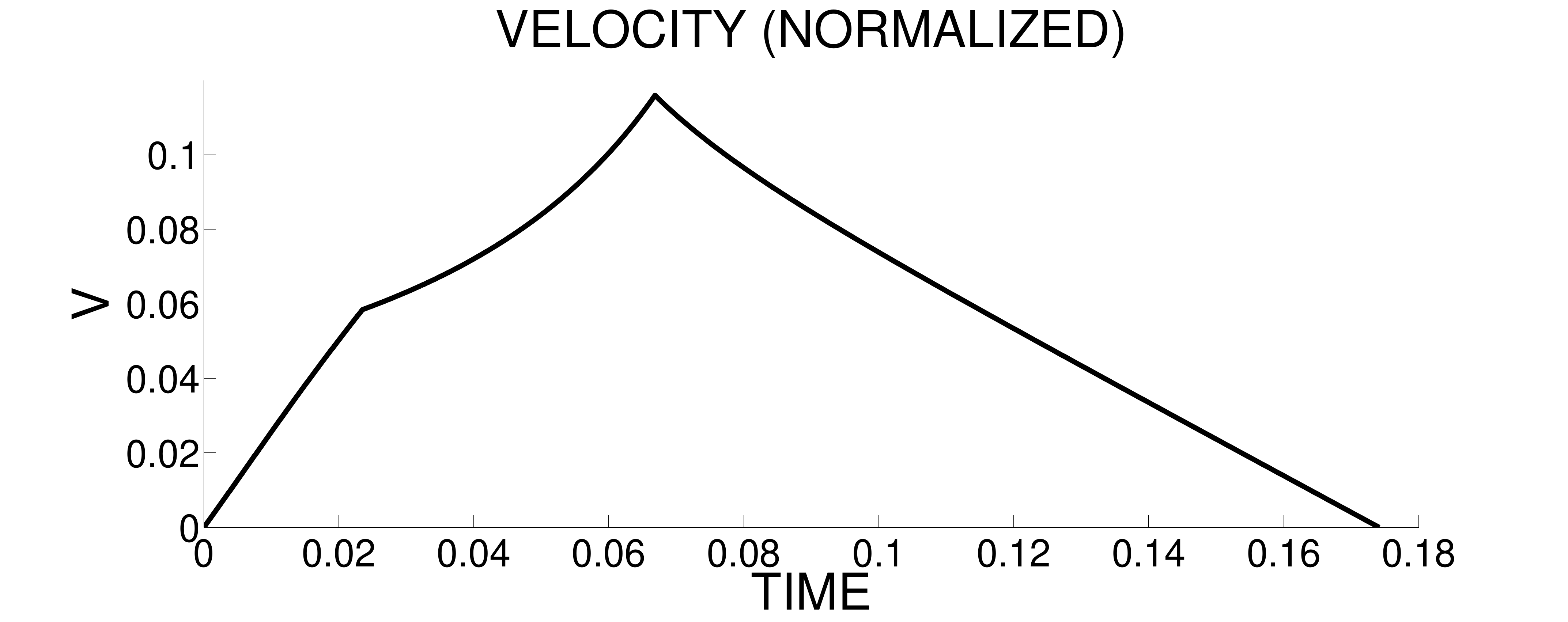}}
\\
\subfigure{\includegraphics [width=5.8cm]{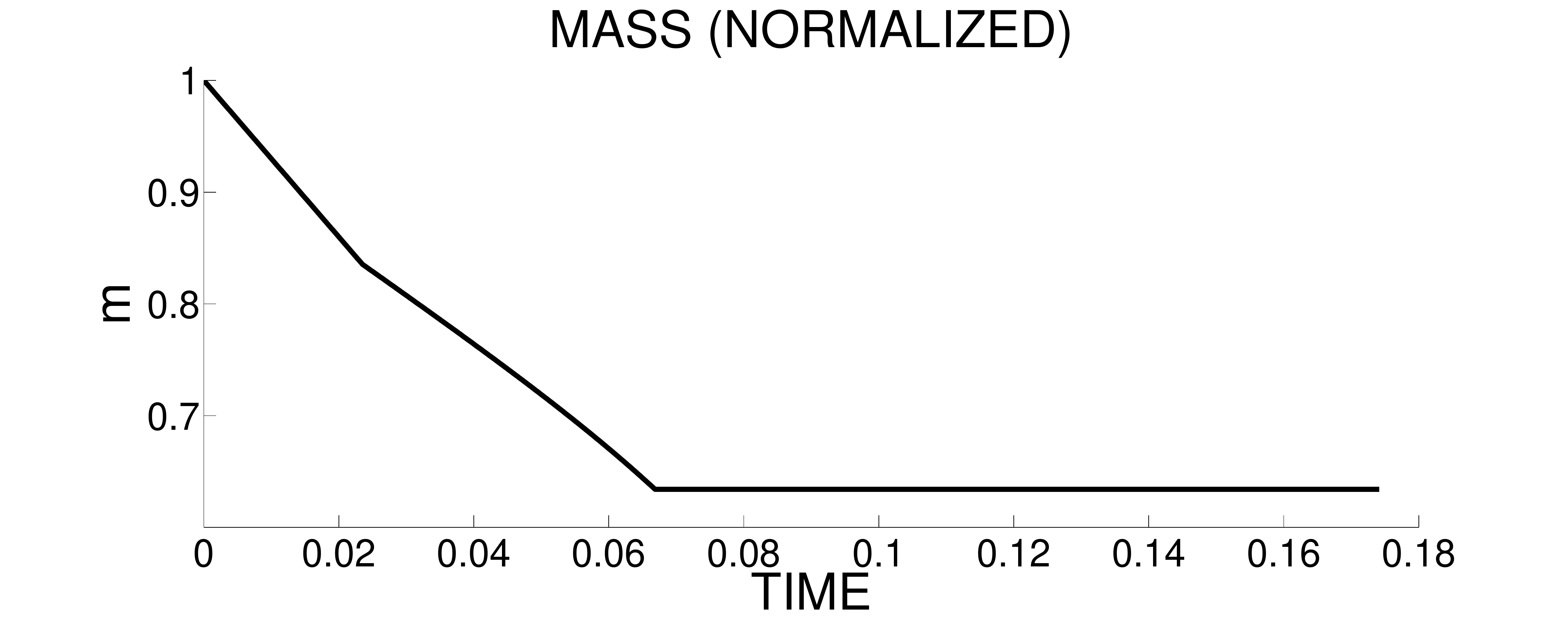}}
\
\subfigure{\includegraphics [width=5.8cm]{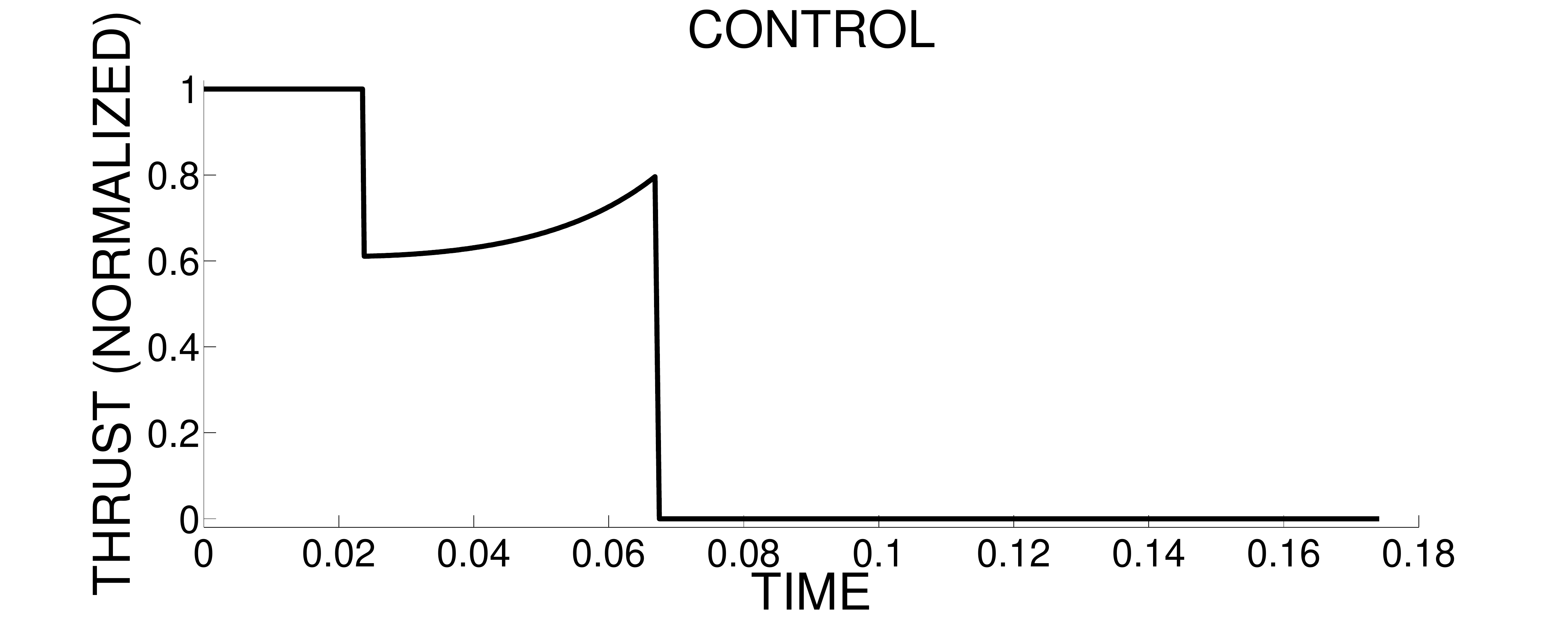}}
\end{center}
\caption{Goddard Problem}
\label{chap2FigGoddard}
\end{figure}

\begin{figure}[!h]
\centering
\begin{center}
\subfigure{\includegraphics [width=5.8cm]{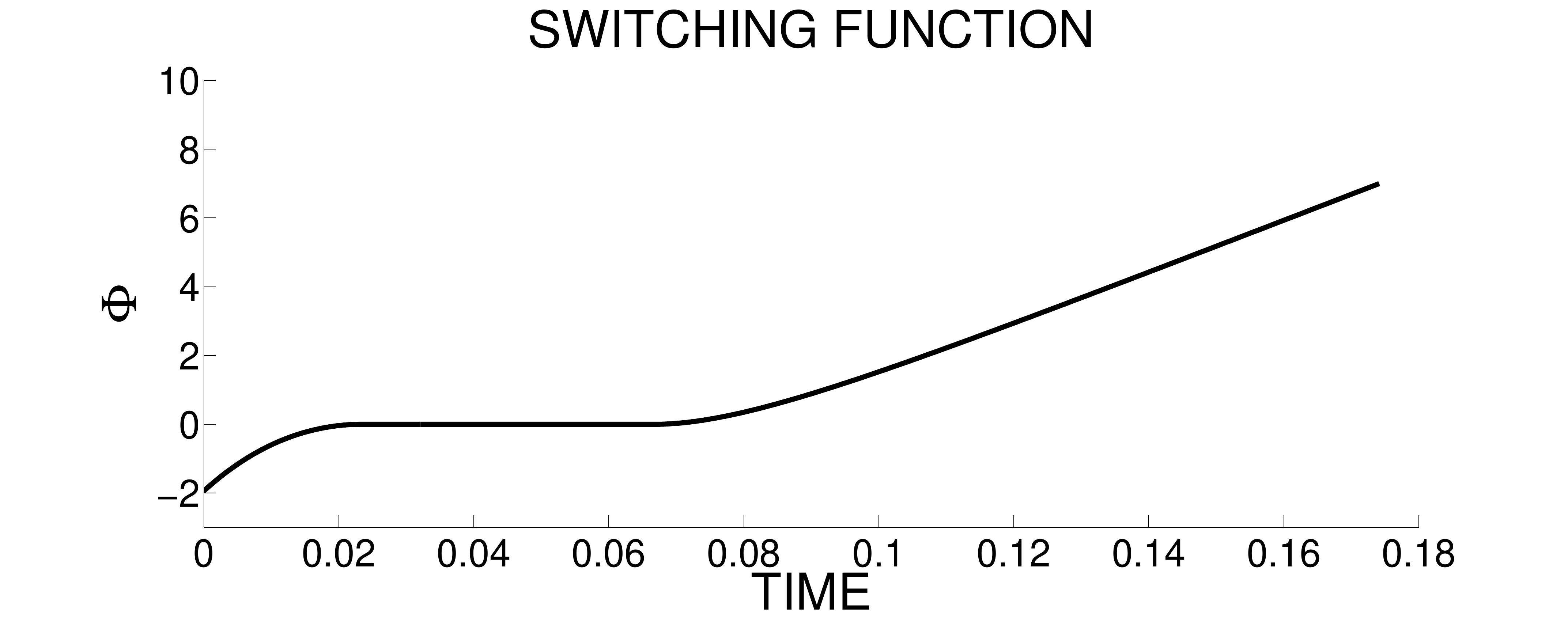}}
\, 
\subfigure{\includegraphics [width=5.8cm]{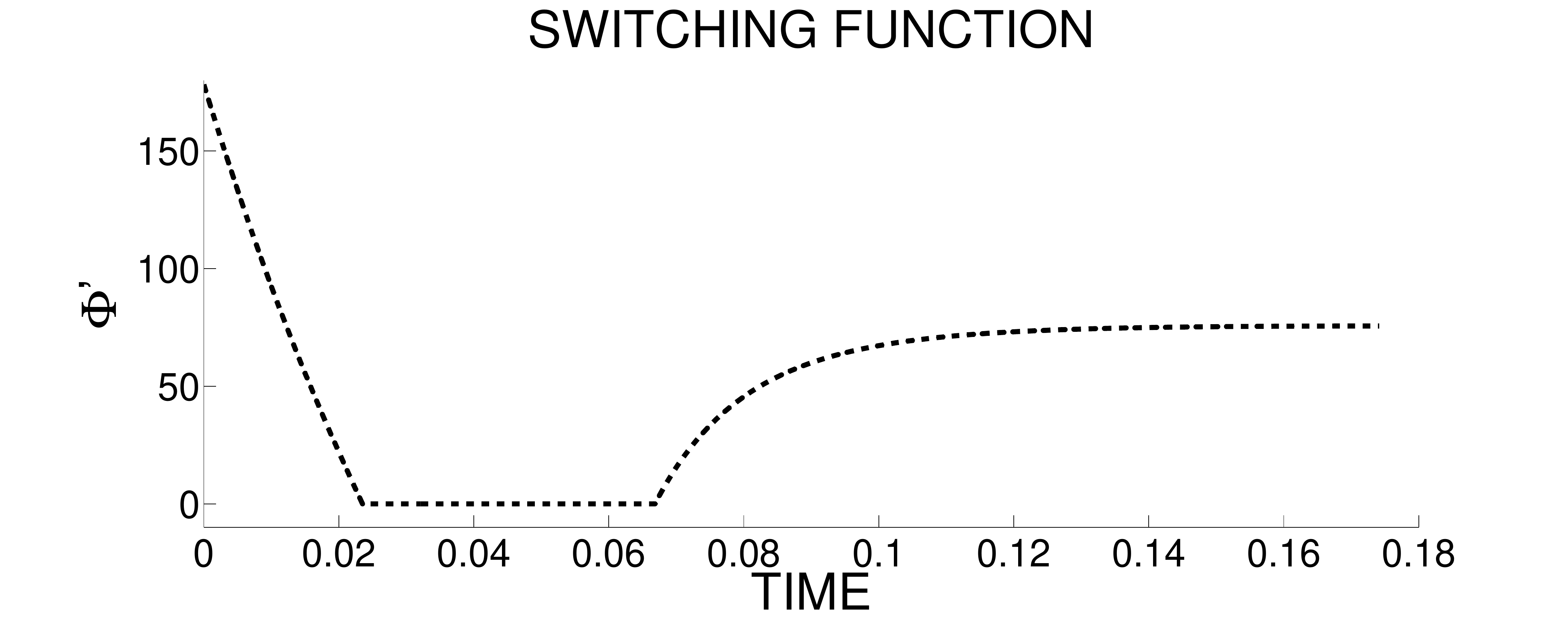}}
\\
\subfigure{\includegraphics [width=7cm]{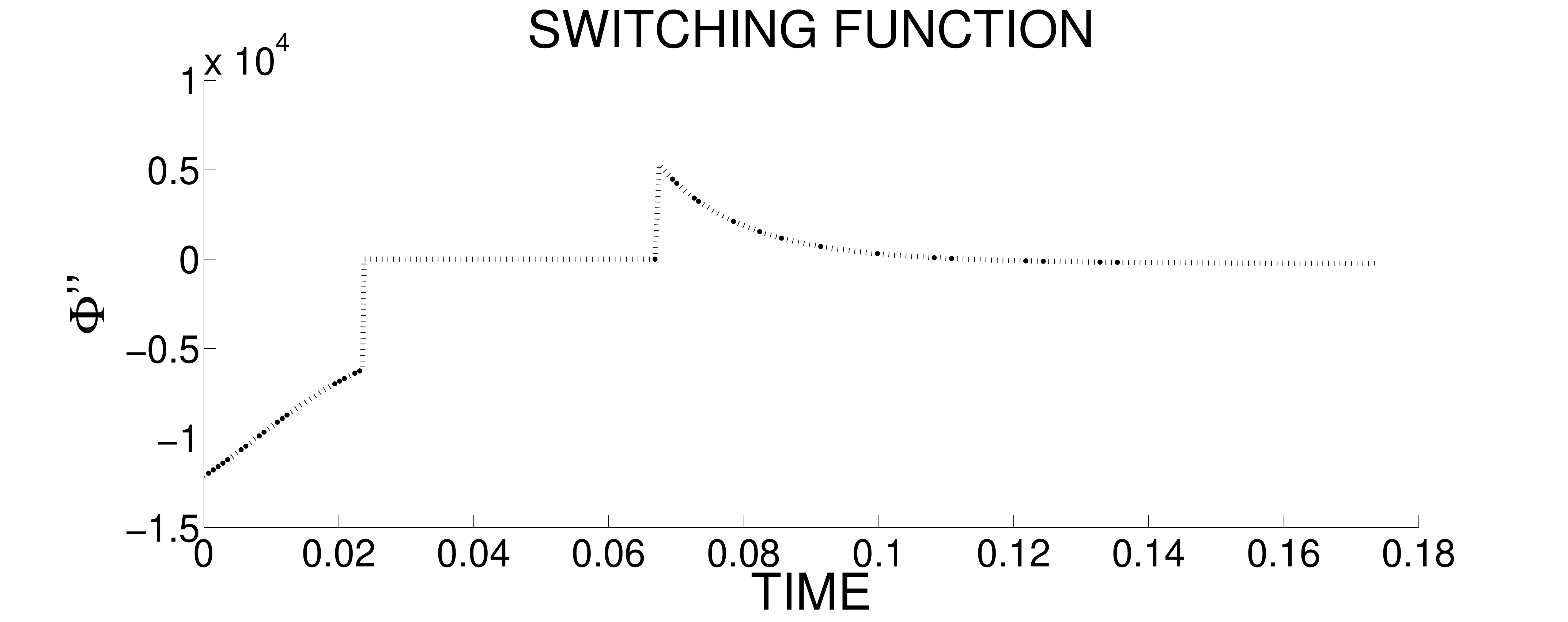}}
\end{center}
\caption{Goddard Problem}
\label{chap2FigGoddard2}
\end{figure}

The third example is the well-known Goddard problem, introduced in Goddard \cite{Goddard} and studied for instance in Seywald-Cliff \cite{SeyCli93}. 
This problem models the ascent of a rocket through the atmosphere, and we restrict here ourselves to  vertical (unidimensional) trajectories.
The state variables are the altitude, speed and mass of the rocket during the flight, for a total dimension of 3.
The rocket is subject to gravity, thrust and drag forces.
The final time is free, and the objective is to reach a certain altitude with a minimal fuel consumption, i.e. a maximal final mass.
\be\tag{P$_3$} 
\left \lbrace
\ba{rl}
&\max\ m_T,\\
&\dot r = v,\\
&\dot v = - {1}/{r^2} + {1}/{m} (\T_{{\rm max}} u - D(r,v)),\\
&\dot m = - b \T_{{\rm max}} u,\\
&0\leq u_t \leq 1,\quad {\rm  a.e.}\ {\rm  on}\ [0,1],\\
&r_0 = 1, \ v_0=0,\ m_0=1,\\
&r_T = 1.01,\quad T\ {\rm  free}, 
\ea
\right .
\ee
with the parameters $b=7$, $\T_{{\rm max}}=3.5$ and the drag given by 
\benl
D(r,v) := 310 v^2 e^{-500(r-1)}.
\eenl
The pre-Hamiltonian function here is
\benl
H:=p_rv+p_v\big[- {1}/{r^2} + {1}/{m} (\T_{{\rm max}} u - D(r,v))\big] - p_mb\T_{\max}u,
\eenl
where $p_r,$ $p_v$ and $p_m$ are the costate variables associated to $r,$ $v$ and $m,$ respectively. The switching function is
\benl
\Phi:= D_uH = \T_{{\rm max}} [-p_mb + {p_v}/{m}].
\eenl
Hence, the bang-bang optimal control is given by
\benl
\left \lbrace
\ba{rl}
\uh_t = 0 \quad {\rm if}\ \Phi_t > 0,\\
\uh_t = 1 \quad {\rm if}\ \Phi_t < 0,
\ea
\right .
\eenl
and the singular control can be obtained by formally solving $\ddot{\Phi}=0$.
The expression of $\uh_{{\rm singular}},$ however, is quite complicated and is not recalled here.
The solution for this problem has the well-known typical structure \textbf{1-singular-0}, as shown on Figures \ref{chap2FigGoddard} and \ref{chap2FigGoddard2}.


\noindent \textbf{Shooting formulations.} Once again fixing the control structure, the shooting unknowns are
\benl
\nu = (p_{1,0}, p_{2,0}, p_{3,0}, t_1, t_2, T) \in \cR^6.
\eenl
Here it is the classical shooting formulation with the entry conditions on $t_1$
\benl
\shoot_3(\nu) := (x_{1,T} - 1.01, p_{2,T}, p_{3,T} + 1, \Phi_{t_1}, \dot \Phi_{t_1}, H_T),
\eenl
while the extended formulation is
\benl
\tilde{\shoot}_3(\nu) := (x_{1,T} - 1.01, p_{2,T}, p_{3,T} + 1, \Phi_{t_1}, \dot \Phi_{t_1}, H_T, [H]_{t_1}, [H]_{t_2}).
\eenl

\subsection{Results} \label{chap2results}

All tests were run on a 12-core platform, with the parallelized (OPENMP) version of the SHOOT (\cite{RR-7380}) package.
The ODE solver is a fixed step 4th. order Runge Kutta method with 500 steps.
The classical shooting is solved with a basic Newton method, and the extended shooting with a basic Gauss-Newton method.
We do not use improved versions of these schemes since we aim to study the behavior of the shooting algorithm in its pure state.
Both algorithms use a fixed step length of 1 and a maximum of 1000 iterations.
In addition to the singular/bang structure, the value of the control on the bang arcs is also fixed according to the expected solution.

The values for the initial costates are taken in $[-10,10]$, and the values for the entry/exit times in $[0,T]$ for $(P_1)$ and $(P_2).$
For $(P_3)$, the entry, exit and final times are taken in $[0,0.2]$.
The number of grid points is set around to 10000 for the three problems.
These grids for the starting points are quite large and rough, which explains the low success rate for $(P_1)$ and $(P_3)$.
However, the solution was found for all three problems.

For each problem, the results are summarized in 3 tables.
The first table indicates the total CPU time for all shootings over the grid, the success rate of convergence to the solution, the norm of the shooting function at the solution, and the objective value.
The second table recalls the solution found by both formulations: initial costate and junction times, as well as final time for $(P_3)$.
The third table gives the singular values for the Jacobian matrix at the solution, as well as its condition number $\kappa := {\sigma_1}/{\sigma_n}$.

We observe that for all three problems $(P_1)$, $(P_2)$ and $(P_3),$ both formulations converge to the same solution, $\hat\nu$ and the objective being identical to more than 6 digits.
The success rate over the grid, total CPU time and norm of the shooting function at the solution are close for both formulations.
Concerning the singular values and condition number of the Jacobian matrix, we note that, for $(P_2),$ the extended formulation has the smallest singular value going from $10^{-8}$ to $1$, thus improving the condition number by a factor $10^8$.
This is caused by the combination of the two entry conditions into a single one that we used in the classical formulation for this problem: as the singular arc lasts until $t_f$, there is only one additional unknown, the entry time.\\

Overall, these results validate the extended shooting formulation, which perform at least as well as the classical formulation and has a theoretical foundation.\\

\begin{remarks} Several additional tests runs were made using the HYBRD (\cite{GaHiMo80}) and NL2SNO (\cite{DenGayWel81}) solvers for the classical and extended shootings instead of the basic Newton and Gauss-Newton method.
The results were similar, apart from a higher success rate for the HYBRD solver compared to NL2SNO.
\end{remarks}

\begin{remarks}  We also tested both formulations using the sign of the switching function to determine the control value over the bang arcs, instead of forcing the value. 
However, this causes a numerical instability at the exit of a singular arc, where the switching function is supposed to be 0 but whose sign determines the control at the beginning of the following bang arc. 
This instability leads to much more erratic results for both shooting formulations, but with the same general tendencies. 
\end{remarks}

\noindent \textbf{Problem 1}\\
Shooting grid: $[-10,10]\times[0,T]^2$, $21^3$ gridpoints, 9261 shootings.

\begin{center}
\begin{tabular}{|l|cccc|}
\hline
Shooting  & CPU & Success & Convergence & Objective\\
\hline
Classical & 74 s  &  21.28 \%  & 1.43E-16  & -106.9059979\\
Extended  & 86 s  &  22.52 \%  & 6.51E-16  & -106.9059979\\
\hline
\end{tabular}\\
\vspace*{0.1cm}
{\bf \small Table 1} $(P_1)$ CPU times, success rate, convergence and objective
\end{center}

\vspace{0.1pt}

\begin{center}
\begin{tabular}{|l|ccc|}
\hline
Shooting  & $p_0$ & $t_1$ & $t_2$\\
\hline
Classical & -0.462254744307241  &  2.37041478456004 &  6.98877992494185 \\ 
Extended  & -0.462254744307242  &  2.37041478456004 &  6.98877992494185 \\
\hline
\end{tabular}\\
\vspace*{0.1cm}
{\bf \small Table 2} $(P_1)$ solution $\hat\nu$ found
\end{center}

\vspace{0.1pt}

\begin{center}
\begin{tabular}{|l|ccc|c|}
\hline
Shooting  & $\sigma_1$ & $\sigma_2$ & $\sigma_3$ & $\kappa$\\ 
\hline
Classical & 3.61 & 0.43 & 5.63E-02 &  64.12\\
Extended  & 27.2 & 1.71 & 3.53E-01 &  77.05\\
\hline
\end{tabular}\\
\vspace*{0.1cm}
{\bf \small Table 3} $(P_1)$ singular values and condition number for the Jacobian
\end{center}

\noindent \textbf{Problem 2} \\
Shooting grid: $[-10,10]^2\times[0,T]$, $21^3$ gridpoints, 9261 shootings.

\begin{center}
\begin{tabular}{|l|cccc|}
\hline
Shooting  & CPU & Success & Convergence & Objective\\
\hline
Classical & 468 s &  94.14 \%   & 1.17E-16 & 0.37699193037\\
Extended  & 419 s &  99.36 \%   & 1.22E-13 & 0.37699193037\\
\hline
\end{tabular}\\
\vspace*{0.1cm}
{\bf \small Table 4} $(P_2)$ CPU times, success rate, convergence and objective
\end{center}

\vspace{0.1pt}
\begin{center}
\begin{tabular}{|l|ccc|}
\hline
Shooting  & $p_{1,0}$ & $p_{2,0}$ & $t_1$\\
\hline
Classical &  0.942173346483640   &     1.44191017584598  &      1.41376408762863     \\ 
Extended  &  0.942173346476773   &     1.44191017581021  &      1.41376408762893 \\
\hline  
\end{tabular}\\
\vspace*{0.1cm}
{\bf \small Table 5} $(P_2)$ solution $\hat\nu$ found
\end{center}

\vspace{1pt}
\begin{center}
\begin{tabular}{|l|ccc|c|}
\hline
Shooting  &  $\sigma_1$ & $\sigma_2$ & $\sigma_3$ & $\kappa$\\ 
\hline
Classical &  24.66 &  5.19 &  1.96E-08  & 1.26E+09 \\
Extended  &  24.70 &  5.97 &  1.13      & 21.86\\
\hline
\end{tabular}\\
\vspace*{0.1cm}
{\bf \small Table 6} $(P_2)$ singular values and condition number for the Jacobian
\end{center}

\noindent \textbf{Problem 3} \\
Shooting grid: $[-10,10]^3\times[0,0.2]^3$, $4^3 \times 5^3$ gridpoints, 8000 shootings.

\begin{center}
\begin{tabular}{|l|cccc|}
\hline
Shooting  & CPU & Success & Convergence & Objective\\
\hline
Classical & 42 s    & 0.82 \%   & 5.27E-13 & -0.634130666 \\
Extended  & 52 s    & 0.85 \%   & 1.29E-10 & -0.634130666 \\
\hline
\end{tabular}\\
\vspace*{0.1cm}
{\bf \small Table 7} $(P_3)$ CPU times, success rate, convergence and objective
\end{center}

\vspace{0.1pt}

\begin{center}
\begin{tabular}{|l|ccc|}
\hline
S.  & $p_{r,0}$ & $p_{v,0}$ & $p_{m,0}$ \\
\hline
C. &    -50.9280055899288    &   -1.94115676279896 &     -0.693270270795148   \\   
E.  &    -50.9280055901093    &   -1.94115676280611 &     -0.693270270787320    \\
\hline
          & $t_1$ & $t_2$ & $t_f$ \\
\hline
C. &    0.02350968417421373 &  0.06684546924474312  & 0.174129456729642     \\ 
E.  &    0.02350968417420884 &  0.06684546924565564  & 0.174129456733106    \\
\hline
\end{tabular}\\
\vspace*{0.1cm}
{\bf \small Table 8} $(P_3)$ solution $\hat\nu$ found (S.: Shooting, C.: Classical, E.: Extended)
\end{center}

\vspace{0.1pt}
\begin{center}
\begin{tabular}{|l|cccccc|c|}
\hline
Shooting  &  $\sigma_1$ & $\sigma_2$ & $\sigma_3$ &  $\sigma_4$ & $\sigma_5$ & $\sigma_6$ & $\kappa$\\ 
\hline
Classical &  6182 &  9.44  &  8.13 &  2.46 &  0.86 &  1.09E-03  & 5.67E+06 \\ 
Extended  &  6189 &  12.30 &  8.23 &  2.49 &  0.86 &  1.09E-03  & 5.67E+06\\ 
\hline
\end{tabular}\\
\vspace*{0.1cm}
{\bf \small Table 9} $(P_3)$ singular values and condition number for the Jacobian
\end{center}


\section{Conclusions}\label{chap2SecConclusion}

Theorems \ref{chap2wp} and \ref{chap2cwp} provide a theoretical support for an extension of the shooting algorithm for problems with all the control variables entering linearly and having singular arcs. The shooting functions here presented are not the ones usually implemented in numerical methods as we have already pointed out in previous section. They come from systems having more equations than unknowns in the general case, while before in practice only square systems have been used. 
Anyway, we are not able to prove the injectivity of the derivative of the shooting function when we remove some equations, i.e. we are not able to determine which equations are redundant, and we suspect that it can vary for different problems. 

The proposed algorithm was tested in three simple problems, where we compared its performance with the classical shooting method for square systems. The percentages of convergence are similar in both approaches, the singular values and condition number of the Jacobian matrix of the shooting function coincide in two problems, and are better for our formulation in one of the problems. Summarizing, we can observe that the proposed method works as well as the one currently used in practice and has a theoretical foundation.

In the  bang-singular-bang case, as in the fishing and Goddard's problems, our formulation coincides with the algorithm proposed by Maurer \cite{Mau76}.

Whenever the system can be reduced to a square one, given that the sufficient condition for the non-singularity of the Jacobian of the shooting function coincides with a sufficient condition for optimality, we could established the stability of the optimal local solution under small perturbations of the data.

\section*{Acknowledgments} We thank the two anonymous referees for their useful remarks.

\end{document}